# Tamed to compatible when $b^{2+} = 1$ and $b^1 = 2$

Clifford Henry Taubes[†]
Department of Mathematics
Harvard University
Cambridge, MA 02138

chtaubes@math.harvard.edu

ABSTRACT: Weiyi Zhang noticed recently a gap in the proof of the main theorem of the author's article *Tamed to compatible: Symplectic forms via moduli space integration* [T] for the case when the symplectic 4-manifold in question has first Betti number 2 (and necessarily self-dual second Betti number 1). This note explains how to fill this gap.

[†] Supported in part by the National Science Foundation

# 1. Introduction

Let X denote a smooth, oriented 4-manifold with the self-dual second Betti number $b^{2+}$ equal to 1. Suppose that X admits a symplectic form giving the orientation; and let ω denote such a form. Theorem 1 in [T] asserted the following:

**Theorem 1**: *Any suitably generic, almost complex structure that is tamed by ω is also compatible with a symplectic form on* X.

The term generic here and in what follows signifies that J can be any almost complex structure from a certain residual set of smooth, ω-tamed almost complex structures. (A residual set is the complement of a countable union of nowhere dense sets. It is, in particular, a dense set in the space of smooth, ω-tamed almost complex structures. Note also that any countable intersection of residual sets is also residual.) By way of a reminder, a complex structure (call it J) is said to be tamed by ω when the bilinear form ω( · ,J( · )) is positive definite. The complex structure is said to be compatible when this bilinear form is also symmetric.

Recently, Weiyi Zhang noted a gap in the proof of this theorem (which was subsequently narrowed by T.-J. Li) in the case when X has first Betti number 2 and when the cup product on $H^1(X;\mathbb{Z})$ is zero. (I am much indebted to Weiyi and T.-J. for their explanations with regards to this matter.) The purpose of this note is to explain why Theorem 1 in [T] holds for this $b^1 = 2$ case also.

# 1. The gap in the proof

The gap involves Proposition 2.1 of [T] which makes a specific assertion about the existence of certain pseudo-holomorphic curves containing some number of chosen points in X. This proposition was claimed in [T] to follow from results in [LL1]-[LL3]. Weiyi Zhang (and T-J. Li) pointed out that the latter results are not sufficient to draw the asserted conclusions when $b^1 = 2$ and the $H^1(X;\mathbb{Z})$ cup product is zero.

**a) Proposition 2.1 in [T]**

The notation needed to state Proposition 2.1 in [T] is as follows: The cup product pairing between classes e and e´ in $H^2(X; \mathbb{R})$ is denoted subsequently by e·e´. A class in $H^2(X; \mathbb{R})$ is said to lie in the *positive cone* when it has positive cup product pairing with itself and with the class that is defined by ω. Introduce c ∈ $H^2(X; \mathbb{Z})$ to denote the first Chern class of the complex line bundle over X whose fiber is the vector subspace in the complexification $\wedge^2 T^*X$ that consists of forms of type (2, 0) as defined by any almost complex structured tamed by ω. (Any two such line bundles are isomorphic.) This class is used to associate the (even) integer $\iota_e = e·e - c·e$ to a given class e ∈ $H^2(X; \mathbb{Z})$.

Let J denote a given ω-tamed almost complex structure. If $C \subset X$ is a J-holomorphic subvariety, then it defines (by integration) an integer valued linear form on $H^2(X; \mathbb{R})$. The Poincaré dual of this class (a class in $H^2(X; \mathbb{Z})$) is denote by $e_C$. Since ω is positive on the smooth part of C, the element $e_C$ is a non-zero class in the positive cone.

Proposition 2.1 in [T] refers to an *irreducible* J-holomorphic subvariety. A subvariety is irreducible if its smooth locus is connected. Any given J-holomorphic subvariety is a union of a finite set of irreducible J-holomorphic subvarieties.

**PROPOSITION 2.1 IN [T]**: *Fix a class* $e \in H^2(X; \mathbb{Z})$ *in the positive cone with* $\hat{e} \cdot \hat{e} > 0$. *There exists* $\kappa > 1 + |c \cdot \hat{e}|/\hat{e} \cdot \hat{e}$ *such that if* n *is an integer greater than* κ, *then the following is true: Suppose that* J *is a given, ω-tamed almost complex structure. Fix* $\frac{1}{2} \iota_{n\hat{e}}$ *points in* X *and there exists a finite set,* Θ, *of pairs of the form* (C, m) *with* $C \subset X$ *an irreducible, J-holomorphic subvariety and* m *a positive integer. Moreover,*
- $\sum_{(C,m) \in \Theta} m\, e_C = n\hat{e}$.
- $\cup_{(C,m) \in \Theta} C$ *contains the chosen set of points*.

To reiterate: W. Zhang pointed out that this proposition has not been established when $b^1 = 2$ and the cup product on $H^1(X; \mathbb{Z})$ is zero.

**b) A replacement for Proposition 2.1 in [T]**

The results in [LL2] and [LL3] do imply a weaker version of the preceding proposition with the assertion being the same except for one less specified point:

**Proposition 1.1**: *Suppose that* X *is a 4-manifold and that* ω *is a symplectic form on* X. *Assume in what follows that the Betti number* $b^{2+}$ *and* $b^1$ *for* X *are* 1 *and* 2 *respectively. Fix a class* $\hat{e} \in H^2(X; \mathbb{Z})$ *in the positive cone with* $\hat{e} \cdot \hat{e} > 0$. *There exists* $\kappa > 1 + |c \cdot \hat{e}|/\hat{e} \cdot \hat{e}$ *such that if* n *is an integer greater than* κ, *then the following is true: Suppose that* J *is a given, ω-tamed almost complex structure. Fix a pair of non-intersecting, embedded curves in* X *that generate* $H_1(X; \mathbb{Z})$/*torsion. Having done this, fix* $\frac{1}{2} \iota_{n\hat{e}}$ - 1 *points in* X *disjoint from the two curves. Given this data, there exists a finite set,* Θ, *of pairs of the form* (C, m) *with* $C \subset X$ *an irreducible, J-holomorphic subvariety and* m *a positive integer. Moreover,*
- $\sum_{(C,m) \in \Theta} m\, e_C = n\hat{e}$.
- $\cup_{(C,m) \in \Theta} C$ *contains the chosen set of* $\frac{1}{2} \iota_{n\hat{e}}$ - 1 *points and at least one point on each of the chosen generators of* $H_1(X; \mathbb{Z})$/*torsion*.

This proposition follows from Lemma 3.3 in [LL2] using the main theorem in [LL3].

## 2. The argument for Theorem 1

Although Proposition 1.1 here is not the same as Proposition 2.1 in [T], it nonetheless leads to Theorem 1 for the case when X has $b^{2+} = 1$ and $b^1 = 2$. (The vanishing or not of the cup product on $H_1(X;\mathbb{Z})$/torsion plays no role in the subsequent arguments.). The argument leading from Proposition 1.1 to Theorem 1 in the case at hand is analogous to what is said in Section 1 of [T]. To set the stage for the argument, let $\Omega$ denote for the moment a given 2-form on X. This form defines a current (a bounded, linear functional on the Frechet space of smooth 2-forms) by the rule

$$\mu \to \Phi_\Omega(\mu) = \int_X \mu \wedge \Omega.$$

(2.1)

If $\Omega$ is a symplectic form and compatible with a given almost complex structure (to be denoted by J), then $\Phi_\Omega$ has the following properties:

- *Fix* $t \in (0,1)$ *and let* $B \subset X$ *denote a ball of radius* t. *Supposing that $\mu$ is a smooth form with support in* B, *then* $|\Phi_\Omega(\mu)| \leq c t^4 \sup_X |\mu|$ *with c being independent of* B *and $\mu$.*
- $\Phi_\Omega$ *annihilates the exterior derivative of 1-forms.*
- $\Phi_\Omega$ *annhilates forms of type (2,0) and (0,2) (as defined by* J).
- *If e is a $\mathbb{C}$-valued type (1,0) form (and* $i = \sqrt{-1}$), *then* $\Phi_\Omega(ie \wedge \bar{e})$ *is bounded from below by a positive number (independent of e) times* $\int_X |e|^2$.

(2.2)

(The first bullet is a consequence only of $\Omega$ being a $C^0$ differential form. The second bullet is a consequence of $\Omega$ being closed. The third bullet follows from the fact that $\Omega$ is a type $(1,1)$ form as defined by J. Given that $\Omega$ is type $(1,1)$, the fourth bullet follows from the fact that $\Omega \wedge \Omega > 0$.)

With (2.2) in mind, the plan in what follows (and this was the plan in [T]) is to approximate the current $\Phi_\Omega$ for a hypothetical symplectic form $\Omega$ by an average of the currents that are defined by integrating 2-forms on (closed) J-holomorphic submanifolds. If a current defined by averaging the integration over closed submanifolds (it is denoted by $\Phi$ below and in Section 1 of [T]) has the properties listed in (2.2), then as explained in Section 1 of [T]), the current $\Phi$ can be smoothed so that it is given by (2.1) with $\Omega$ being a symplectic form compatible with J.

With regards to the conditions in (2.2): If C is a closed, 2-dimensional submanifold in X, then integration over C defines a closed current (it obeys the second bullet in (2.2)); and if C is J-holomorphic, then this current annihlates $(2,0)$ and $(0,2)$ forms (the third bullet in (2.2)). Moreover, integration over a J-holomorphic submanifold of a 2-form such as $ie \wedge \bar{e}$ with $e$ being type $(1,0)$ is always non-negative, which is part of the constraint in the fourth bullet of (2.2). The whole of the fourth bullet is obtained by averaging the currents that are defined by a suitable family of J-holomorphic submanifolds. In particular, the family must be large enough so that the following condition is met:

*Every line in some open cone in the (1,0) tangent space at any given point in X is the J-complex tangent line of some subvariety in the family.*

(2.3)

The first bullet in (2.2) is also a constraint on the averaging family of submanifolds. This bullet says (roughly) that the J-complex tangent lines to the submanifolds in the family can not accumulate to a great extent in any given complex direction at any given point.

The plan just outlined talked about *averaging* the integration currents that are defined by family of J-holomorphic subvarieties. To define this averaging, the family of J-holomorphic submanifolds is chosen (in part) to have a readily available measure so that averaging means integration with respect to this measure. For this purpose (and for other reasons), the family will be parametrized in part by a finite set of distinct points in X. This finite set is denoted by $\Xi$. The family of submanifolds parametrized by a given $x \in \Xi$ is denoted by $\mathcal{K}_x$. This family $\mathcal{K}_x$ is defined with help of a J-holomorphic submanifold containing x (to be denoted by $C_x$) whose fundamental class in $H_2(X; \mathbb{Z})$ is dual to some large n version of the class $n\hat{e}$. (Note in this regard that the purpose of Proposition 2.1 is to obtain (assuming J is generic) a suitable $C_x$ for any $x \in X$.) The submanifolds in $\mathcal{K}_x$ are parametrized by a data set that is written as $D_{x2} \times \mathcal{O}_{x2} \times \mathcal{B}_x$ with

- $D_{x2}$ *is a small radius disk about* 0 *in the normal bundle to* $C_x$ *at* x.
- $\mathcal{O}_{x2}$ *is a small radius disk about the complex line in* $\mathbb{P}(T_{1,0}X)$ *at* x *defined by* $T_{1,0}C$.
- $\mathcal{B}_x$ *is the product of small radius balls about* $\frac{1}{2} \iota_{n\hat{e}}$ -2 *distinct points in* $C_x$−x.

(2.4)

The parametrization map is denoted by $\Gamma_x$; it is a map from $D_{x4} \times \mathcal{O}_{x4} \times \mathcal{B}_x$ onto an open neighborhood of $C_x$ in the space of J-holomorphic submanifolds of X.

The desired current $\Phi$ is defined as follows: Supposing that $\upsilon$ is a smooth 2-form,

$$\Phi(\upsilon) = \sum_{x \in \Xi} \int_{\{q \in D_{x2} \times \mathcal{O}_{x2} \times \mathcal{B}_x\}} \left( \int_{\Gamma_x(q)} \upsilon \right)$$

(2.5)

Note that the averaging family that is used here is different than that used in [T]. As it turns out, this family is simpler to use; and it could have been used in [T] to simplify some of the latter's constructions.

Where does the genericity assumption on J arise? It is used in what follows to find a suitable $C_x$ for each $x \in X$. The submanifold $C_x$ must be generic with regards to certain maps from spaces of J-holomorphic submanifolds with marked points; and if J is itself generic, then an appropriately generic $C_x$ can be found for each x.

The remaining sections of this note explain how to choose, for each $x \in X$, the required submanifold $C_x$ with the disks $D_{x2}$ and $\mathcal{O}_{x2}$ and the ball $\mathcal{B}_x$, and how to define the associated parametrizing map $\Gamma_x$. These sections also explain the criteria for choosing the finite subset $\Xi$. Given this data, the upcoming Proposition 5.1 asserts in effect that (2.5)'s current obeys the conditions in (2.2) (only the first and fourth bullets need to be checked.) As noted already, the almost complex structure J is tamed by a symplectic form when (2.2) holds.

## 3. Moduli spaces and transversality

This section explains how to choose an appropriate J-holomorphic submanifold $C_x$ for each $x \in X$. This is where Proposition 1.1 is invoked; and it is where the genericity conditions on J are used.

### a) Notation

Some additional notions from [T] are needed for what follows. To start the introduction, suppose that $e \in H^2(X;\mathbb{Z})$ is any given class. Write $\iota_e$ (which is e·e - c·e) as 2d; and introduce g to denote $\frac{1}{2}(e·e + c·e) + 1$. The significance of these numbers is as follows. With regards to g: If J is a given almost complex structure and if C is an irreducible, J-holomorphic *submanifold* in X with $e_C = e$, then g is the genus of C. If C is irreducible but with singularities, then the genus of its model curve (this is defined in Section 2b in [T?]) is strictly less than g. (This model curve genus is denoted here by $k_C$.) With regards to the integer d (it is always an integer): Given $k \in \{0, 1, \ldots, g\}$, let $\mathcal{M}_{e,k}$ denote the space of irreducible, J-holomorphic subvarieties $C \subset X$ with $e_C = e$ and with $k_C = k$. In the case $k_C = g$, every element in $\mathcal{M}_{e,g}$ is a J-holomorphic submanifold and the number 2d is the index of the operator whose kernel is the vector space of first order deformations of C as a J-holomorphic submanifold. This is the operator $D_C$ that appears in (2.12) of [T]. If J is generic, then $\mathcal{M}_{e,g}$ is a smooth manifold of dimension 2d and its tangent space at C is canonically isomorphic to the kernel of $D_C$ (the cokernel is 0 then).

Letting e again denote a given class in $H^2(X;\mathbb{Z})$, Proposition 3.1 in [T] considers a space it denotes by $\mathcal{M}_e$ which is defined as follows: Any given element in $\mathcal{M}_e$ (call it $\Theta$) is a finite set of pairs where each pair has the form (C,m) with C being an irreducible, J-holomorphic subvariety in X and with m being a positive integer. The set $\Theta$ is constrained by the condition that

$$\sum_{(C,m) \in \Theta} m \, e_C = e.$$
(3.1)

As noted by Proposition 3.1 in [T], the space $\mathcal{M}_e$ is compact with respect to the topology of pointwise convergence as subsets of X and convergence as currents. (This is nothing more than Gromov's compactness theorem from [G].) In this regard, the current defined by $\Theta$ is given by the rule:

$$\upsilon \to \Sigma_{(C,m)\in\Theta} \, m \int_C \upsilon \ .$$

(3.2)

The space $\mathcal{M}_{e,g}$ embeds in $\mathcal{M}_e$ as an *open* subspace via the map sending any given submanifold C to the set consisting of the single pair (C, 1). Given a class e, and a set $\Theta \in \mathcal{M}_e$, use [$\Theta$] to denote the set $\{(e_C, k_C)\}_{(C,m)\in\Theta}$. Then, use $\mathfrak{N}_e$ to denote the collection of such sets [$\Theta$]. The set $\mathfrak{N}_e$ is always finite (because $\mathcal{M}_e$ is compact).

Here is a convention that is used henceforth: The propositions that are stated in this section and subsequently assume without saying that a given class $\hat{e} \in H^2(X; \mathbb{Z})$ is chosen from the positive cone with $\hat{e}\cdot\hat{e} > 0$. The class e in what follows is taken to be some very large integer multiple of $\hat{e}$. The numbers d and g used here are defined using this choice for e (thus $d = \frac{1}{2}(e\cdot e - c\cdot e)$ and $g = \frac{1}{2}(e\cdot e + c\cdot e) + 1$.)

**b) Submanifolds in $\mathcal{M}_{e,g}$ and points in X**

The proposition that follows is the first step to finding the desired submanifold $C_x$ for each point in x in X.

**Proposition 3.1**: *Suppose that J is generic. Each point in X is contained in some submanifold from $\mathcal{M}_{e,g}$.*

The upcoming proof of Proposition 3.1 and subsequent proofs require various space of J-holomorphic subvarieties in X with marked points and/or with other constraints. What follows directly introduces those that are relevant for the proof of Proposition 3.1.

To start this introduction, fix an element [$\Theta$] is from the set $\mathfrak{N}_e$ and introduce the associated space $\mathcal{M}_{[\Theta]} = \times_{(e',k')\in[\Theta]} \mathcal{M}_{e',k'}$. To be sure, an element in $\mathcal{M}_{[\Theta]}$ is a finite set of irreducible, J-holomorphic subvarieties. The union of these is therefore a reducible J-holomorphic subvariety. An element in $\mathcal{M}_{[\Theta]}$ and the corresponding union (in X) of the constituent subvarieties (thus, a reducible subvariety in X) will not be notationally distinguished in what follows.

Propositions 3.3 and 3.4 in [T] have the following implication: If J is generic, then $\mathcal{M}_{[\Theta]}$ is a smooth manifold whose dimension is at most 2d - 2 except in the event that [$\Theta$] is the one element set $\{(e, g)\}$ (in which case $\mathcal{M}_{[\Theta]}$ is $\mathcal{M}_{e,g}$ and its dimension is 2d). It is assumed henceforth that J is chosen so that all such $\mathcal{M}_{[\Theta]}$ spaces are smooth manifolds.

A second general notion that is needed for what follows is that of an *image variety*. This is defined in Section 3.1 of [T] as follows: Let Y denote a manifold and let n denote a positive integer no greater than the dimension of Y. A subset $\mathcal{S} \subset Y$ is a dimension n image variety (or a codimension dim(Y) - n image variety) if it is closed, and if each point in Y has a neighborhood that intersects $\mathcal{S}$ as the image via a smooth, proper map of a manifold with a finite number of connected components, none with dimension greater than n (or codimension less than dim(Y) - n).

*Proof of Proposition 3.1*: An embedded circle in X will be called a *loop*. Fix a loop $\gamma$ and then, having fixed an element $[\Theta]$ from $\mathfrak{N}_e$, let $\mathcal{M}^\gamma_{[\Theta]}$ denote the subspace in $(\times_{(e',k') \in [\Theta]} \mathcal{M}_{e',k'}) \times X$ consisting of elements that have the form (C, t) with $t \in \gamma \cap C$. With $\gamma$ fixed, arguments much like those in the Sections 5.1-5.4 of [T] (Sections Ac-Ad in [T*]) prove that this space is image variety inside $(\times_{(C,m) \in \Theta} \mathcal{M}_{e_C, k_C}) \times X$. In this case, none of the manifolds that define this image variety $\mathcal{M}^\gamma_{[\Theta]}$ have dimension greater than 2d - 3 unless $[\Theta]$ is the one element set $\{(e, g)\}$ in which case $\mathcal{M}^\gamma_{e,g}$ is a smooth manifold whose dimension is 2d -1.

By the same token, if N is a positive integer, and if $\{\gamma_1, \ldots, \gamma_N\}$ is a finite set of pairwise disjoint loops, then each $[\Theta] \in \mathfrak{N}_e$ and $i \in \{1, \ldots, N\}$ version of $\mathcal{M}^{\gamma_i}_{[\Theta]}$ will be a smooth manifold if J is generic. (But note that the loops have to be specified before J.)

The preceding can be generalized as follows: Let r denote a positive integer, and let $\mathcal{X}_r$ denote the open set of points in $\times_r X$ with pairwise distinct entries. Having chosen a loop $\gamma$ and an element $[\Theta] \in \mathfrak{N}_e$, let $\mathcal{M}^\gamma_{[\Theta],r}$ denote the subset of $\mathcal{M}^\gamma_{[\Theta]} \times \mathcal{X}_r$ consisting of elements that have the form ((C, t), $\mathfrak{v}$) with each entry of $\mathfrak{v}$ being a point in C−t. If J is generic, then this too is a smooth image variety. (As before, the loop $\gamma$ determines what is meant by 'generic'.) The maximum dimension of the relevant manifolds for this image variety is 2r more than the maximal dimension of those that define $\mathcal{M}^\gamma_{[\Theta]}$. (All of this can be proved using the techniques in Sections 5.1-5.4 in [T] which is Sections Ac-Ae in [T*].) In the case when $[\Theta]$ is the one point set $\{(e,g)\}$, then the corresponding space $\mathcal{M}^\gamma_{e,g,r}$ is a smooth manifold of dimension 2d - 1 + 2r. In any case, let $\pi_r$ denote the map from $\mathcal{M}^\gamma_{[\Theta],r}$ to $\mathcal{X}_r$ that is induced by the projection from $\mathcal{M}^\gamma_{\Theta} \times \mathcal{X}_r$.

If N is a positive integer, and if a set of N loops is specified (call this set $\Gamma$), then each $\gamma \in \Gamma$ version of $\mathcal{M}^\gamma_{[\Theta],r}$ will be an image variety when J is generic. Moreover, assuming J is generic, then the corresponding N versions of the maps $\pi_r$ will be in general position with respect to each other. This is to say the following: Label the loops in $\Gamma$ as $\{\gamma_1, \ldots, \gamma_N\}$ and introduce $\mathcal{M}^\Gamma_{[\Theta]}$ to denote the subset in $(\times_{1 \leq i \leq N} \mathcal{M}^{\gamma_i}_{[\Theta]}) \times \mathcal{X}_r$ consisting of elements $((C_1, t_1), \ldots, (C_N, t_N), \mathfrak{v})$ with each $(C_i, t_i)$ in the corresponding $\mathcal{M}^{\gamma_i}_{[\Theta]}$, with $\mathfrak{v}$ from $\mathcal{X}_r$, and with each entry of $\mathfrak{v}$ contained in each $C_i$ and distinct from the point $t_i$. This space $\mathcal{M}^\Gamma_{[\Theta]}$ is itself an image variety (when J is generic) whose constituent manifolds have dimension no greater than the expected dimension from naïve dimension counting. (This can be proved using the same techniques as those in the Appendix of [T] or [T*].)

Of particular interest with regards to Proposition 3.1 is the case when $r = d-1$ and $[\Theta]$ is not $\{(e,g)\}$. In this case, any loop $\gamma$ version of $\mathcal{M}^\gamma_{[\Theta],d-1}$ is an image variety of dimension $4d-5$ or less; and then, given N disjoint loops, the corresponding $\mathcal{M}^\Gamma_{[\Theta]}$ is an image variety of dimension $4d-4-N$ or less. (This is what is expected from dimension counting for the intersection of N submanifolds of codimension 1 or more in an ambient $4d-4$ dimensional space.) The point now is that $\mathcal{M}^\Gamma_{[\Theta],d-1}$ must be empty when N is greater than $4d-4$ (because its dimension would be negative in this case). What this means when $N > 5d-2$ is this: Given any $d-1$ tuple of distinct points in X (call this $d-1$ tuple $\mathfrak{v}$), there is at least one loop (call it $\gamma$) from the set $\Gamma$ such that no element from $\mathcal{M}_{[\Theta]}$ that intersects $\gamma$ contains all entries of $\mathfrak{v}$. By taking N larger still, there will be at least one loop (called $\gamma$ again) from the set $\Gamma$ such that no element from $\mathcal{M}_{[\Theta]}$ that intersects $\gamma$ contains all entries of $\mathfrak{v}$ when $[\Theta]$ is *any* element from $\mathfrak{N}_e$ except possibly the element $\{(e,g)\}$.

Keeping the preceding in mind, suppose now that all of the loops in $\Gamma$ represent the same generator of $H_1(X;\mathbb{Z})/$torsion. If all of the entries of the given $\mathfrak{v}$ are not in any one subvariety from some $\gamma \in \Gamma$ version of $\mathcal{M}^\gamma_{[\Theta]}$ when $[\Theta]$ is not $\{(e,g)\}$, then it follows from Proposition 1.1 that all of the entries of $\mathfrak{v}$ must be in a submanifold from $\mathcal{M}^\gamma_{(e,g)}$.

To complete the proof of Proposition 3.1, fix $x \in X$ and then fix another $d-2$ distinct points in X that are distinct from x and from each loop in $\Gamma$. Set $\mathfrak{v} = (x, \mathfrak{w})$. Then by virtue of what is said in the preceding paragraph, there exists $C \subset \mathcal{M}_{e,g}$ containing x.

**b) The map $\pi_{d-1}$**

Supposing that r is a positive integer, let $\mathcal{M}_{e,g,r}$ denote the subspace of pairs in $\mathcal{M}_{e,g} \times \mathcal{X}_r$ that have form $(C, \mathfrak{v})$ with each entry of $\mathfrak{v}$ being a point in C. If J is generic (as will be assumed), then this is a smooth submanifold in $\mathcal{M}_{e,g} \times \mathcal{X}_r$ of dimension $2(d+r)$ with the projection to $\mathcal{M}_{e,g}$ being a fibration. The latter map is denoted by $\pi_\mathcal{M}$. Meanwhile, the map to $\mathcal{X}_r$ from $\mathcal{M}_{e,g,r}$ induced by the projection from $\mathcal{M}_{e,g} \times \mathcal{X}_r$ is denoted by $\pi_r$. A point $(C, \mathfrak{v})$ in $\mathcal{M}_{e,g,r}$ is a regular value of $\pi_r$ when its differential at $(C, \mathfrak{v})$ has trivial cokernel. By way of a parenthetical remark, this notion of regular value for $\pi_r$ has the following reinterpretation: Introduce the operator $D_C$ from (2.12) in [T] or [T*]. Since J is generic, this operator has kernel dimension $2d$ and cokernel dimension 0. Let $\ker(D_{C,\mathfrak{v}})$ denote the vector subspace in $\ker(D_C)$ consisting of the elements that vanish at each entry of $\mathfrak{v}$. This vector subspace has dimension at least $2(d-r)$. The element $(C, \mathfrak{v})$ is a regular point of $\pi_r$ if and only if this dimension is exactly $2(d-r)$.

Of particular interest are the cases when $r = 1$ and when $r = d-1$.

**Proposition 3.2**: *Suppose that J is generic. Given x, there is a submanifold $C \in \mathcal{M}_{e,g}$ containing x such that the pair (C,x) is a regular point of the map $\pi_1$ on $\mathcal{M}_{e,g,1}$. More to the point, there exists $\mathfrak{w} \in \mathcal{X}_{d-2}$ with each entry distinct from x and contained in C, and such that the pair $(C, (x, \mathfrak{w}))$ (which is in $\mathcal{M}_{e,g,d-1}$) is a regular point of the map $\pi_{d-1}$.*

It is important to note that if $(C, (x, \mathfrak{w}))$ is a regular point of $\pi_{d-1}$ in $\mathcal{M}_{e,g,d-1}$, then (C, x) is a regular point of $\pi_1$ in $\mathcal{M}_{e,g,1}$ and $(C, \mathfrak{w})$ is a regular point of $\pi_{d-2}$ in $\mathcal{M}_{e,g,d-2}$.

Keep in mind for what follows that the condition that $(C, (x, \mathfrak{w}))$ be a regular point of $\pi_{d-1}$ is an open condition on the set of points $\mathcal{M}_{e,g,d-1}$: There is an open neighborhood of $(C, (x, \mathfrak{w}))$ in $\mathcal{M}_{e,g,d-1}$ consisting entirely of regular points for $\pi_{d-1}$. What is more, this is an open condition with regards to small deformations of the almost complex structure J.

*Proof of Proposition 3.2*: The proof has two parts. The first part explains why (when J is generic) there are submanifolds in $\mathcal{M}_{e,g}$ containing x which partner with x to make for a regular point in $\mathcal{M}_{e,g,1}$ of the map $\pi_1$. The second part proves the rest of Proposition 3.2.

*Part 1*: To start: Techniques from the appendix to [T] or [T*] can be readily used to prove the following: If J is generic, then the set of critical points of $\pi_r$ in $\mathcal{M}_{e,g,r}$ is an image variety whose constituent manifolds have dimension at most 4r - 1. (This is the expected dimension: The evaluation map from $\ker(D_C)$ to $\oplus_{w: w \text{ is an entry of } \mathfrak{w}} N_C|_w$ can be viewed as a $2d \times 2r$ matrix; and the subspace of $2d \times 2r$ dimensional matrices with non-zero cokernel is a stratified space whose strata have codimension no less than $2(d-r)+1$.) In particular, if J is generic, then the set of critical points of $\pi_1$ in $\mathcal{M}_{e,g,1}$ has dimension 3. Denote this set by $C_1$. If x is a regular value of $\pi_1$, then x is not in $\pi_1(C_1)$ so if C is in $\mathcal{M}_{e,g}$ and C contains x, then (C,x) is a regular point of $\pi_1$.

Now suppose that x is critical value of $\pi_1$ (so it is in $\pi_1(C_1)$). To deal with this case, suppose for the moment that $(C, \mathfrak{v})$ is an element in $\mathcal{M}_{e,g,d-1}$. Write $\mathfrak{v}$ as $(p, \mathfrak{w})$ with p being a point in X and $\mathfrak{w}$ being in $\mathcal{X}_{d-2}$. Mapping $(C, \mathfrak{v})$ to (C, p) defines a map from $\mathcal{M}_{e,g,d-1}$ to $\mathcal{M}_{e,g,1}$. This map is denoted by $f_1$. The space $f_1^{-1}(C_1)$ is an image variety in $\mathcal{M}_{e,g,d-1}$ whose constituent manifolds have dimension at most $2d+1$ (which is 2d - 2 more than 3). Let $\pi_{d-2}$ denote (also) the map from $\mathcal{M}_{e,g,d-1}$ that sends any given $(C, (p, \mathfrak{w}))$ to $\mathfrak{w}$. The $\pi_{d-2}$ image of $f_1^{-1}(C_1)$ in $\mathcal{X}_{d-2}$ is also a $2d+1$ dimensional image variety. Since $\dim(\mathcal{X}_{d-2}) = 4d - 8$ which is greater than $2d+1$ when d is large, there are points in $\mathcal{X}_{d-2}$ from the complement of $\pi_{d-2}(f^{-1}(C_1))$ (when d is large) with no entry equal to x. This has the following implication: Let $\mathfrak{w}$ denote such a point and let C denote a submanifold from $\mathcal{M}_{e,g}$ containing x and all entries of $\mathfrak{w}$. (The proof of Proposition 3.1 explains why C exists.) The element $(C, (x, \mathfrak{w}))$ is in $\mathcal{M}_{e,g,d-1}$ but since it isn't in $f^{-1}(C_1)$, the pair (C,x) is a regular value of $\pi_1$.

*Part 2*: A simple dimension counting argument such as that used in Part 1 won't work to prove the second part of Proposition 3.2 because critical locus of $\pi_{d-1}$ on $\mathcal{M}_{e,g,d-1}$ is a codimension 3 image variety when J is generic, which is to say that its maximal dimension is $4d-5$. Since this is larger than $\dim(\mathcal{X}_{d-2})$, it might be that every point in $\mathcal{X}_{d-2}$ is in the $\pi_{d-2}$ image of the $\pi_{d-1}$ critical locus.

To prove the second assertion of Proposition 3.2, first fix submanifold $C \subset \mathcal{M}_{e,g}$ that contains x and is such that $(C,x)$ is a regular value of $\pi_1$. Remember: This means that $\ker(D_{C,x})$ has dimension $2d-2$. (The space $\ker(D_{C,x})$ is the vector subspace of $\ker(D_C)$ that vanish at x.) Fix a linearly independent pair $\{s_1, s_2\}$ from $\ker(D_{C,x})$. Then, there is an open, dense set of points in C where this pair spans $N_C$. (To prove this, suppose to the contrary that there is an open set in C where $s_1 = rs_2$ with r being a real valued function. By virtue of $s_1$ and $s_2$ obeying (2.12) in [T] or [T*], this function r obeys $\bar{\partial} r = 0$. Since r is real, it must obey $dr = 0$; and so it is constant. This implies that $s_1 - rs_2$ is zero on an open set, which makes it everywhere zero since elements in the kernel of $D_C$ vanish at exactly e·e points (counting multiplicities). The vanishing of $s_1 - rs_2$ is nonsense because $s_1$ and $s_2$ are chosen to be linearly independent elements of $\ker(D_C)$.). With the preceding understood, fix a point $w_1 \in C$ where $\{s_1, s_2\}$ span $N_C$. By virtue of this, the kernel of the evaluation map $s \to s(w_1)$ on $\ker(D_{C,x})$ has dimension $2d-4$. Now chose a linearly independent pair $\{s_3, s_4\}$ from $\ker(D_{C,(x,w_1)})$. Arguing as with $\{s_1, s_2\}$, there exists $w_2 \in C$ where $\{s_3, s_4\}$ span the normal bundle to C. This implies that the evaluation map $s \to s(w_2)$ on $\ker(D_{C,(x,w_1)})$ has kernel dimension $2d-6$, so this is the dimension of $\ker(D_{C,(x,w_1,w_2)})$. Continuing in this vein constructs a $(d-2)$-tuple $\mathfrak{w} = (w_1, \ldots w_{d-2})$ of points in C such that the $\ker(D_{C,)x,\mathfrak{w}})$ has dimension 2. This implies that $(C,(x,\mathfrak{w}))$ is a regular point of $\pi_{d-1}$.

### c) Maps to $\mathbb{P}(T_{1,0}X)$

Let $\mathbb{P}(T_{1,0}X)$ denote the $\mathbb{CP}^1$ fiber bundle over X whose fiber at any given point is the space of complex lines in the $(1,0)$ tangent space at the point. There is a canonical section over $\mathcal{M}_{e,g,1}$ of the $\pi_1$ pull-back of this bundle: It assigns to any given pair $(C,p)$ the $(1,0)$ tangent line to C at the point p. This section is denoted in what follows by $\psi$. The bundle $\pi_1^*(\mathbb{P}(T_{1,0}X)$ and the section $\psi$ play a central role in this section.

Supposing that $\mathfrak{w}$ is a point in $\mathcal{X}_{d-2}$, let $\mathcal{M}^{*\mathfrak{w}}$ denote the subspace of regular points of $\pi_{d-2}$ on $\mathcal{M}_{e,g,d-2}$ that map to $\mathfrak{w}$. This is the space of pairs $(C, \mathfrak{w})$ with each entry of $\mathfrak{w}$ in C and with $\dim(\ker(D_{C,\mathfrak{w}})) = 4$. It is a smooth, 4-dimensional manifold; and, since $\mathfrak{w}$, is fixed it can (and will) be considered a submanifold in $\mathcal{M}_{e,g}$. Let w now denote an entry of

w. Then the section $\psi$ defines a smooth map (to be denoted by $\psi^w$) from $\mathcal{M}^{*w}$ to $\mathbb{P}(T_{1,0}X)|_w$ that sends any given submanifold C from $\mathcal{M}^{*w}$ to its (1,0) tangent line at $w$.

**Proposition 3.3**: *Suppose that J is generic. Given* $x \in X$, *there exists* $w \in X_{d-2}$ *with entries distinct from X and a submanifold C from* $\mathcal{M}_{e,g}$ *containing x and each entry of w, and such that*:
- *The triple* (C, (x,w)) *is a regular point of* $\pi_{d-1}$ *on* $\mathcal{M}_{e,g,d-1}$.
- *If w is an entry of w, then C (which is in* $\mathcal{M}^{*w}$) *is a regular point of* $\psi^w$ *on* $\mathcal{M}^{*w}$.

Note that these conditions are also open conditions on the data (C,(x,w)) (and also on the choice of almost complex structure.)

*Proof of Proposition 3.3*: After an appeal to Proposition 3.1, choose (C´, w´) from $\mathcal{M}_{e,g,d-2}$ so that (C´,(x,w´)) is a regular point of $\pi_{d-1}$. Any (C,w) from $\mathcal{M}_{e,g,d-2}$ near this initial choice such that $x \in C$ will have (C,(x,w)) being a regular point of $\pi_{d-1}$. This said, then the first bullet of the proposition will be satisfied if (C,w) obeys the second, it is close to (C´,w´), and x is in C.

Keeping the preceding in mind, let $\mathcal{M}^*_{e,g,d-2}$ denote the subspace of regular values of $\pi_{d-2}$ in $\mathcal{M}_{e,g,d-2}$. If J is generic, then this is the complement of a dimension 5 image variety. Define a subspace $\mathcal{F}_{d-2} \subset \mathcal{M}^*_{e,g,d-2}$ by the following rule: A pair (C,w) is in $\mathcal{F}_{d-2}$ if there is an entry of w (call it w) such that C is a critical point in $\mathcal{M}^{*w}$ of the map $\psi^w$. Note in this regard that $C \in \mathcal{M}^{*w}$ is a critical point of $\psi^w$ if and only if the linear map sending $\ker(D_{C,w})$ to $(N_C \otimes T^{1,0}C)|_w$ by the rule $s \to \partial s|_w$ is not surjective (so its kernel dimension is at least 3). Here, $\partial s$ is the (1,0) derivative of s. The techniques in the appendix to [T] or [T*] can be used to prove that $\mathcal{F}_{d-2}$ is a codimension 3 image variety if J is generic. (This is the expected dimension because the subspace of 2×4 matrices with kernel dimension 3 or more is a stratified space with top stratum having codimension 3.) Codimension 3 means that the constituent manifolds of $\mathcal{F}_{d-2}$ have dimension at most 4d - 7 (because $\dim(\mathcal{M}_{e,g,d-2})$ is 4d - 4).

Now, given $x \in X$, let $\mathcal{M}^{*x}_{d-2}$ denote the subspace of $\mathcal{M}^*_{e,g,d-2}$ that consists of the pairs (C, w) with C containing x and with x not an entry of w. If (C,w) is in $\mathcal{M}^{*x}_{d-2}$ and (C,(x,w)) (which is in $\mathcal{M}_{e,g,d-1}$) is a regular point for $\pi_{d-1}$, then $\mathcal{M}^{*x}_{d-2}$ near (C,w) is a smooth 4d - 6 dimensional submanifold (this is 2 less than the dimension of $\mathcal{M}^*_{d-2}$). Since this dimension is one greater than that of $\mathcal{F}_{d-2}$, there are points (C,w) in $\mathcal{M}^x_{d-2}$ in any given neighborhood of the chosen pair (C´,w´) that are not in $\mathcal{F}_{d-2}$. Any such (C,w) which is close to (C´,w´) obeys the conditions of Proposition 3.3.

Suppose now that x ∈ X and that (C,𝔴) with x are described by the two bullets in Proposition 3.3. The first bullet implies that the subset of submanifolds in $\mathcal{M}^{*\mathfrak{w}}$ near C that contain x is a 2-dimensional submanifold of $\mathcal{M}^{*\mathfrak{w}}$ near C. To put this in a broader context, define a subspace $\mathcal{M}^{*(x,\mathfrak{w})}$ in $\mathcal{M}^{*\mathfrak{w}}$ as follows: A submanifold (call it C´) from $\mathcal{M}^{*\mathfrak{w}}$ is in this space when it contain x and when the evaluation map from $\ker(D_{C',\mathfrak{w}})$ to $N_C|_x$ is surjective. This is a smooth, 2-dimensional submanifold of $\mathcal{M}^{*\mathfrak{w}}$ that contains the submanifold C, and it contains all submanifolds in $\mathcal{M}^{*\mathfrak{w}}$ near C that contain x.

Since all of the constituent submanifolds in $\mathcal{M}^{*(x,\mathfrak{w})}$ contain x, the section ψ restricts to define a map from $\mathcal{M}^{*(x,\mathfrak{w})}$ to $\mathbb{P}(T_{1,0}X)|_x$. This map is denoted by $\psi^x$.

**Proposition 3.4**: *Suppose that J is generic. Given* x ∈ X, *there exists* 𝔴 ∈ $X_{d-2}$ *with entries distinct from* X *and a submanifold C from* $\mathcal{M}_{e,g}$ *containing* x *and each entry of* 𝔴, *and such that*:
- *The triple* (C, (x,𝔴)) *is a regular point of* $\pi_{d-1}$ *on* $\mathcal{M}_{e,g,d-1}$.
- *If w is an entry of* 𝔴, *then C (which is in* $\mathcal{M}^{*\mathfrak{w}}$) *is a regular point of* $\psi^w$ *on* $\mathcal{M}^{*\mathfrak{w}}$.
- *The submanifold C (which is in* $\mathcal{M}^{*(x,\mathfrak{w})}$) *is a regular value of* $\psi^x$ *on* $\mathcal{M}^{*(x,\mathfrak{w})}$.

As was the case with Proposition 3.3, all of these conditions are open conditions on the data (C,(x,𝔴)) (and also on the choice of almost complex structure.)

***Proof of Proposition 3.4***: Choose a pair C´ ∈ $\mathcal{M}_{e,g}$ and 𝔴´ ∈ $X_{d-2}$ with entries disjoint from x and contained in C´ with the requirement that (C´, (x,𝔴´)) obey the conditions stated by the first two bullets of Proposition 3.3. Any pair (C,𝔴) near to (C´,𝔴´) in $\mathcal{M}_{e,g,d-2}$ with x ∈ C will also obey the two bullets of Proposition 3.3. Since these are the first two bullets of Proposition 3.4, the task is to find such a (C, 𝔴) near to (C´,𝔴´) that obeys the third bullet of the proposition. Keep in mind here that any given C ∈ $\mathcal{M}^{*(x,\mathfrak{w})}$ is a regular point of the restriction of $\psi^x$ to $\mathcal{M}^{*(x,\mathfrak{w})}$ when the map from $\ker(D_{C,(x,\mathfrak{w})})$ to $N_C|_x$ sending s to ∂s is surjective.

To find the desired pair (C,𝔴), consider first the subspace $\mathcal{M}^*_{e,g,1} \subset \mathcal{M}_{e,g,1}$ consisting of the regular points for the map $\pi_1$. Let $D\pi_1$ denote the differential of $\pi_1$. The kernel of $D\pi_1$ is the subbundle in $T\mathcal{M}^*_{e,g,1}$ whose fiber at any given (C,p) is $\ker(D_{C,p})$ (these are the sections that vanishes at p). The differential of ψ along the fibers of $\pi_1$ defines a bundle homomorphism over $\mathcal{M}^*_{e,g,1}$ from $\ker(D\pi)$ to the 2-dimensional vector bundle whose fiber at (C, p) is $N_C|_p$ (the normal bundle to C at p). This differential is denoted by Dψ. It maps a 2d - 2 dimensional vector bundle to a 2 dimensional vector bundle. Let $\mathcal{F}$ now denote the subspace in $\mathcal{M}^*_{e,g,1}$ where Dψ is not surjective. The techniques from the appendix of [T] or [T*] can be used to prove the following: If J is generic, then $\mathcal{F}$ is a 5-dimensional image variety in $\mathcal{M}^*_{e,g,1}$. Assume henceforth that J is generic so that dim($\mathcal{F}$) ≤ 5.

Now, given the point x, the subspace $\pi_1^{-1}(x)$ in $\mathcal{M}^*_{e,g,1}$ is a manifold of dimension 2d-2 which is much greater than 5 when d is large. Therefore, there exists submanifolds in $\mathcal{M}_{e,g}$ as close as desired to C´ that contain x with (C,x) $\notin \mathcal{F}$. Choose such a submanifold C with this property that is very close to C´. By virtue of (C,x) being in $\mathcal{M}^*_{e,g,1}$, the space ker($D_{C,x}$) has dimension 2d-2. By virtue of (C,x) not being in $\mathcal{F}$, there is a pair of linearly independent sections $\{s_1, s_2\} \subset \ker(D_{C,x})$ with $\{\partial s_1, \partial s_2\}$ spanning $N_C|_x$. Let $K \subset \ker(D_{C,x})$ denote the kernel of the map $s \to \partial s|_p$. This is a 2d-4 dimensional subspace. Fix a linearly independent pair $\{s_3, s_4\}$ in K. Mimick the argument in Part 2 of the proof of Proposition 3.2 to find a point $w_1$ in C very near to the first entry of $\mathfrak{w}´$ (assuming C is very near to C´) where $\{s_3, s_4\}$ spans the normal bundle to C. Having chosen $w_1$, let $K_1 \subset K$ denote the kernel of the evaluation map $s \to s(w_1)$. This space has dimension 2d-6. Choose linearly independent $\{s_5, s_6\}$ from $K_1$ and a point $w_2$ in C very near to the second entry of $\mathfrak{w}´$ (assuming again that C is very close to C´) where this pair spans $N_C$. Continue in this vein (just like in Part 2 of the proof of Proposition 3.2) to generate a point $\mathfrak{w} \in X_{d-2}$ very near to $\mathfrak{w}´$ (assuming as before that C is very close to C´). The resulting pair (C, $\mathfrak{w}$) with x obeys all three bullets of Proposition 3.4.

**d) Maps to X**

Given $x \in X$, Proposition 3.4 supplies pairs (C, $\mathfrak{w}$) from $\mathcal{M}_{e,g,d-2}$ with various nice properties. But, there is yet one last requirement for the (C, $\mathfrak{w}$); this is stated in the fourth bullet of the next proposition.

**Proposition 3.5**: *Suppose that J is generic. Given $x \in X$, there exists $\mathfrak{w} \in X_{d-2}$ with entries distinct from X and a submanifold C from $\mathcal{M}_{e,g}$ containing x and each entry of $\mathfrak{w}$, and such that*:
- *The triple (C, (x,$\mathfrak{w}$)) is a regular point of $\pi_{d-1}$ on $\mathcal{M}_{e,g,d-1}$.*
- *If w is an entry of $\mathfrak{w}$, then C (which is in $\mathcal{M}^{*\mathfrak{w}}$) is a regular point of $\psi^w$ on $\mathcal{M}^{*\mathfrak{w}}$.*
- *The submanifold C (which is in $\mathcal{M}^{*(x,\mathfrak{w})}$) is a regular value of $\psi^x$ on $\mathcal{M}^{*(x,\mathfrak{w})}$.*
- *If p is any point from C (but not an entry of $\mathfrak{w}$), then the evaluation map from ker($D_{C,\mathfrak{w}}$) to $N_C|_p$ is a surjection.*

Noting that dim(ker($D_{C,\mathfrak{w}}$)) = 4 and dim($N_C|_p$) = 2, one expects that in a generic situation, the conditions in the fourth bullet of Proposition (3.5) can be met for any x. This is indeed the case: Items e) and f) of the second bullet of Proposition 3.7 in [T] or [T*] imply that a pair (C, $\mathfrak{w}$) described by Proposition 3.5 can be found for any given x. The argument for this is outlined momentarily. What follows directly is an observation about the evaluation map form ker($D_{C,\mathfrak{w}}$) to $N_C|_p$ for p near x and also for p near an entry of $\mathfrak{w}$.

The first point to make is that the fourth bullet of Proposition 3.5 holds when p is near x and also when p is near any given entry of $\mathfrak{w}$ when $(C, \mathfrak{w})$ obeys the first three bullets (which just restate Proposition 3.4). Indeed, the fourth bullet of Proposition 3.5 holds for p near x by virtue of the first bullet. As explained directly, the fourth bullet also holds for p near any given entry of $\mathfrak{w}$ by virtue of the second bullet. To see this is, let $w$ denote an entry of $\mathfrak{w}$. Complex coordinates (call then $(z,t)$) can be found for a neighborhoood of $w$ in X so that $w$ is the origin, so that the $z$ = constant slices are J-holomorphic disks, and so that the $t = 0$ slice is part of C. This implies, first, that the span of $\frac{\partial}{\partial t}$ at $t = 0$ can be identified with the normal bundle to C. With this identification understood, then any given section from $\ker(D_{C,\mathfrak{w}})$, call it $\eta$, appears in these coordinates as a function of $z$ and $\bar{z}$ having the form

$$\eta(z) = \alpha z + \mathcal{O}(|z|^2) \quad and \quad \bar{\partial}\eta = \alpha + \mathcal{O}(|z|)$$

(3.3)

with $\alpha \in \mathbb{C}$. By virtue of the second bullet in Proposition 3.5, there exists sections $\eta_1$ and $\eta_i$ from $\ker(D_{C,\mathfrak{w}})$ such that the corresponding versions of $\alpha$ are the complex numbers 1 and i. Now, if p is a point from C near $w$, then p has $(z, t)$ coordinates $z = z$ and $t = 0$ with z being some small normed complex number. Then $\eta_1$ at p is $z(1 + \mathcal{O}(|z|))$ and $\eta_i$ at p is $iz(1+ \mathcal{O}(|z|))$; and so if $|z|$ is very small, then these span the normal bundle to C at p.

The preceding observation about the a priori surjectivity of the evaluation map on $\ker(D_{C,\mathfrak{w}})$ near entries of $\mathfrak{w}$ has the following implication:

**Lemma 3.6**: *If $(C,(x,\mathfrak{w})) \in \mathcal{M}_{e,g,d\text{-}1}$ obeys the conclusions of Proposition 3.5, then they are also obeyed by all points on some open neighborhood in $\mathcal{M}_{e,g,d\text{-}1}$ of $(C,x,\mathfrak{w})$.*

(3.4)

*Proof of Lemma 3.6*: As noted previously, the conditions in the first three bullets of Proposition 3.5 are open conditions. The condition in the fourth bullet is not a priori an open condition because the differential of $\pi_1$ on $T\mathcal{M}^{*\mathfrak{w}}|_C$ at a point (C, p) is zero at the entries of $\mathfrak{w}$. However, given small radius balls in X centered on the entries of $\mathfrak{w}$, then there is an open neighborhood of $(C,(x,\mathfrak{w}))$ in $\mathcal{M}_{e,g,d\text{-}1}$ consisting of data $(C´, (x´,\mathfrak{w}´))$ where the fourth bullet of Proposition 3.5 holds using $(C´,(x´, \mathfrak{w}´))$ in lieu of $(C,(x,\mathfrak{w}))$ when p (which is supposed to be a point in C´) is not in any of the chosen balls. Meanwhile, as noted subsequent to (3.3), if the second bullet of Proposition 3.5 holds without regard for the fourth bullet, then there is a punctured ball in X centered at each entry of $\mathfrak{w}$ and a neighborhood of $(C,(x,\mathfrak{w}))$ in $\mathcal{M}_{e,g,d\text{-}1}$ consisting of data $(C´, (x´,\mathfrak{w}´))$ where the fourth bullet of Proposition 3.5 holds using $(C´,(x´, \mathfrak{w}´))$ in lieu of $(C,(x,\mathfrak{w}))$ when p (which is supposed to be a point in C´) is in any one of these punctured balls.

What follows directly is a four part outline of the proof that the fourth bullet condition of Proposition can be satisfied. The genericity assertions made below are established in the Appendix to [T] or [T*].

*Part 1*: To set notation, supposing that k is a positive integer, let $\mathcal{M}^*_{e,g,k} \subset \mathcal{M}_{e,g,k}$ denote the set of regular points for $\pi_k$. Keep in mind that any $(C, (w_1, \ldots, w_k))$ from $\mathcal{M}_{e,g,k}$ is a critical point of $\pi_k$ when the evaluation map from $\ker(D_C)$ to $\oplus_k N_C|_{w_k}$ is not surjective. As noted previously, if J is generic, then the set of critical points of $\pi_k$ is a codimension $2(d-k)+1$ image variety, hence a variety with maximal dimension at most $4k-1$. Therefore, the set of critical *values* of $\pi_k$ in $\mathcal{X}_k$ is also an image variety with top dimension at most $4k-1$. The case $k = d - 2$ is of initial interest. In this case, the set of critical points of $\pi_{d-2}$ is an image variety of dimension $4d-9$ when J is generic. As a consequence, so is the set of critical values of $\pi_{d-2}$. (This number $4d - 9$ is one less than the dimension of $\mathcal{X}_{d-2}$.) The take-away here is that $\pi_{d-2}^{-1}(\mathfrak{w})$ is in $\mathcal{M}^*_{e,g,d-2}$ for a residual set of elements $\mathfrak{w}$.

*Part 2*: The subspace $\mathcal{M}^*_{e,g,d-2}$ is introduced in part because the dimension of the kernel of the differential of $\pi_{d-2}$ on $\mathcal{M}^*_{e,g,d-2}$ is constant. More to the point, this kernel defines a 4-dimensional vector bundle over $\mathcal{M}^*_{e,g,d-2}$ to be denoted by $\ker(D\pi_{d-2})$.

Let $\mathcal{M}^{**}_{e,g,d-1} \subset \mathcal{M}_{e,g,d-1}$ denote the set of elements that have the form $(C,(p,\mathfrak{w}))$ with $(C,\mathfrak{w})$ in $\mathcal{M}^*_{e,g,d-2}$. This is an open set in $\mathcal{M}_{e,g,d-1}$. The space $\mathcal{M}^{**}_{e,g,d-1}$ has a subspace (it is denoted by $\mathcal{L}$) that is defined as follows: Let $\pi_1$ denote the map sending $(C,(p,\mathfrak{w}))$ to p. An element $(C,(p,\mathfrak{w}))$ is in $\mathcal{L}$ when the differential of $\pi_1$ on $\ker(D\pi_{d-2})$ is not surjective (so its kernel dimension is 3 or 4).

If J is generic, then $\mathcal{L}$ is a codimension 3 image variety. Thus, its strata have dimension at most $4d-5$ because $\mathcal{M}^{**}_{e,g,d-1}$ has dimension at most $4d-2$. Meanwhile, the set of regular values of the restriction of $\pi_{d-2}$ to each constituent manifold of $\mathcal{L}$ is the complement of a measure zero set. As a consequence, there is a residual set in $\mathcal{X}_{d-2}$ consisting of regular values of $\pi_{d-2}$ on each constituent manifold of $\mathcal{L}$. Here is what this means: Supposing that $\mathfrak{w}$ is a point in $\mathcal{X}_{d-2}$, let $\mathcal{M}^{*\mathfrak{w}}_1$ denote the subset of pairs $(C,p)$ with C in $\mathcal{M}^{*\mathfrak{w}}$ and with p in C (but not an entry of $\mathfrak{w}$). (Note that this is space is denoted by $\mathcal{M}^{*\mathfrak{w}}_X$ in [T] or [T*].) Now suppose that $\mathfrak{w}$ is a regular value of $\pi_{d-2}$ on the constituent manifolds of $\mathcal{L}$. Then, the space of pairs $(C,p) \in \mathcal{M}^{*\mathfrak{w}}_1$ with $(C,(p, \mathfrak{w}))$ in $\mathcal{L}$ is also an image variety in $\mathcal{M}^{*\mathfrak{w}}_1$ with codimension at 3 or more, which is dimension 3 or less. (The dimension of $\mathcal{M}^{*\mathfrak{w}}_1$ is six.)

With the preceding understood, given x, no generality is lost by assuming that the pair $(C, \mathfrak{w})$ that appears in Proposition 3.4 has its $\mathfrak{w}$ being a regular value of $\pi_{d-2}$ on $\mathcal{L}$. (Remember that the conditions in Proposition 3.4 are open conditions on the data, so if the initial choice of $\mathfrak{w}$ is not a regular value of $\pi_{d-2}$ on $\mathcal{L}$, then any sufficiently generic nearby choice for $\mathfrak{w}$ will be.)

*Part 3*: Let $\mathcal{L}^*$ denote the subset of $\mathcal{M}^*_{e,g,d}$ consisting of elements $(C,(x,p,\mathfrak{w}))$ such that $(C,(x,\mathfrak{w}))$ (which is in $\mathcal{M}_{e,g,d-1}$) is a regular value of $\pi_{d-1}$ and such that $(C,(p,\mathfrak{w}))$ is in $\mathcal{L}$. This space $\mathcal{L}^*$ is the inverse image of $\mathcal{L}$ via the forgetful map, so it is a codimension 3 image variety in $\mathcal{M}_{e,g,d}$. (The forgetful map in this case sends $(C,(x,p,\mathfrak{w}))$ to $(C,(p,\mathfrak{w}))$.) Therefore, the strata of $\mathcal{L}^*$ have dimension at most $4d-3$ because $\mathcal{M}_{e,g,d}$ has dimension $4d$.

Let $\mathcal{L}^{*\lozenge}$ denote the (open) subset of $\mathcal{L}^*$ whose elements have the form $(C,(x,p,\mathfrak{w}))$ with p different from x and from each entry of $\mathfrak{w}$. This set maps to X by the rule $(C,(x,p,\mathfrak{w})) \to x$. Denote this map from $\mathcal{L}^{*\lozenge}$ to X by $\pi_\lozenge$. For generic J, the set of critical points of $\pi_\lozenge$ will be an image variety with codimension no less than $4d-6$, which means that its constituent manifolds have dimension at most 3. To elaborate, note that the tangent space to a constituent manifold of $\mathcal{L}^{*\lozenge}$ containing $(C, (x, p, \mathfrak{w}))$ consists of data that has the form $(\eta, v_x, \cdots)$ where $\eta \in \ker(D_C)$, and $v_x$ is in $TX|_x$, and where the terms indicated by $\cdots$ are not relevant to what is going to be said except to the extent that they and $\eta$ obey linear relations that do not depend on $v_x$. These relations are written collectively as $S(\eta, \cdots) = 0$. Meanwhile, $\eta$ and $v_x$ obey the linear relation $\eta|_x = \Pi_C v_x$ with $\Pi_C$ denoting the orthogonal projection to the normal bundle of C at x. This is the only constraint on $v_x$. The differential of $\pi_\lozenge$ sends $(\eta, v_x, \cdots)$ to $v_x$. As a consequence, the differential is surjective if and only if the set of $(\eta, \cdots)$ obeying the $S(\eta, \cdots) = 0$ maps surjectively to $N_C|_x$ via the map $(\eta, \cdots) \to \eta(x)$. This is a map from a vector space of dimension $4d-5$ to a vector space of dimension 2. The non-surjective linear maps from $\mathbb{R}^{4d-5}$ to $\mathbb{R}^2$ form a real analytic subspace whose top stratum has codimension $4d-6$. If J is generic, then this is also the minimum codimension of the set of critical points of $\pi_\lozenge$ on $\mathcal{L}^{*\lozenge}$. (See Part 13 of Section 5.6 in [T] which is Section Af in [T*] for a full proof of this dimension 3 assertion.)

*Part 4*: If J is generic, then it follows from what is said in the preceding paragraph that the image inside $\mathcal{X}_{d-2}$ of the critical point set of $\pi_\lozenge$ on $\mathcal{L}^{*\lozenge}$ via the map sending $(C,(x,p,\mathfrak{w}))$ to $\mathfrak{w}$ is itself an image variety of dimension at most 3. This image variety is denoted in what follows by $\mathfrak{Z}$. In particular, if $\mathfrak{w}$ is not in this image variety $\mathfrak{Z}$, then any $(C,(x,p,\mathfrak{w}))$ from $\mathcal{L}^{*\lozenge}$ is a regular point for $\pi_\lozenge$. With this understood, then it follows as a consequence that $(C,\mathfrak{w})$ in Proposition 3.4 can be chosen so that $\mathfrak{w}$ is not in this set $\mathfrak{Z}$.

To continue, if $(C,(x,p,\mathfrak{w}))$ is a regular point of $\pi_\diamond$ (for example, if $\mathfrak{w}$ is not in $\mathfrak{Z}$), then $\pi_\diamond^{-1}(x)$ near $(C,(x,p,\mathfrak{w}))$ (which is a subset of $\mathcal{L}^{*\diamond}$) is an image variety whose constituent manifolds have dimension at most $4d-7$. This implies that the image of $\mathcal{L}^{*\diamond}$ in $\mathcal{M}^*_{e,g,d-2}$ near $(C,\mathfrak{w})$ via the map that sends $(C,(x,p,\mathfrak{w}))$ to $(C,\mathfrak{w})$ is also a $4d-7$ dimensional image variety. With this number $4d-7$ in mind, let $\mathcal{M}^x_{e,g,d-2}$ denote the subset of $\mathcal{M}_{e,g,d-2}$ consisting of pairs $(C,\mathfrak{w})$ with C containing x. This space has dimension $4d-6$. As a consequence, there is a dense set in $\mathcal{M}^x_{e,g,d-2}$ consisting of pairs $(C,\mathfrak{w})$ that aren't from $\mathcal{L}^{*\diamond}$. And, if $(C,\mathfrak{w})$ is not from $\mathcal{L}^{*\diamond}$, then (by definition) $\ker(D\pi_{d-2})$ generates the fiber of $N_C$ at each point on C except at the entries of $\mathfrak{w}$.

## 4. The definition of the set $\mathcal{K}$

The definition of the set $\mathcal{K}$ has six steps which are presented in Sections 4a-4f.

### a) Step 1 of $\mathcal{K}$'s definition

What follows is the first step to defining $\mathcal{K}$:

STEP 1: *Given* $x \in X$, *choose* $\mathfrak{w}_x \in X_{d-2}$ *and* $C_x \in \mathcal{M}_{e,g}$ *so that the bullets in Proposition 3.5 hold when* $\mathfrak{w} = \mathfrak{w}_x$ *and* $C = C_x$.

It proves convenient to define a complex coordinate system on a neighborhood of the point x using the submanifold $C_x$. This neighborhood is denoted by $B_x$; it is identified by the coordinate chart with a product of disks (a polydisk). The polydisk in this case is written as $\Delta_x \times D_x$ with $\Delta_x$ being a disk in $C_x$ centered at x and with $D_x$ being a disk about the zero point in the normal bundle to $C_x$ at x. The coordinates on $\Delta_x \times D_x$ are denoted by $(z, t)$ and the embedding of $\Delta_x \times D_x$ to X is denoted by $\lambda_x$. If the radius of $\Delta_x$ and $D_x$ are sufficiently small (as will henceforth be assumed), then $\lambda_x$ can be described as follows:

- *The point* $(0, 0)$ *is mapped by* $\lambda_x$ *to the point* x.
- *The disk* $\Delta_x \times \{0\}$ *is mapped by* $\lambda_x$ *holomorphically onto an open set of* x *in* $C_x$.
- *Each* $z \in \Delta_x$ *version of the disk* $\{z\} \times D_x$ *is mapped by* $\lambda_x$ *to a J-holomorphic disk*.
- *The pull-back by* $\lambda_x$ *of* $T^{1,0}X$ *is spanned at each point by* $\mathbb{C}$-*valued 1-forms*

$$\sigma^0 = dz + \mathfrak{a}\, d\bar{z} \quad and \quad \sigma^1 = dt + \mathfrak{b}\, d\bar{t} + \mathfrak{c}\, d\bar{z}$$

*where* $\mathfrak{a}, \mathfrak{b}$ *and* $\mathfrak{c}$ *are* $\mathbb{C}$-*valued functions that vanish at* $t = 0$.

(4.1)

Here is a schematic picture:

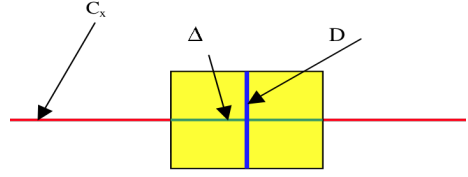

The part of $C_x$ in green is $\Delta_x$. The blue line is $D_x$. The yellow region is the image of $\lambda_x$

### b) Step 2 of $\mathcal{K}$'s definition

To set things up for Step 2, reintroduce the $\mathbb{CP}^1$-bundle $\mathbb{P}(T_{1,0}X)$ over X. Let $\mathfrak{p}$ denote the projection map. It is convenient in what follows to view $\mathbb{P}(T_{1,0}X)$ as the space of pairs (p, L) with $p \in X$ and with L being a complex line in $T_{1,0}X|_p$. This view facilitates the definition of a map (to be denoted by $\psi^{\mathfrak{w}_x}$) from $\mathcal{M}_1^{*\mathfrak{w}_x}$ to $\mathbb{P}(T_{1,0}X)$ that covers the projection induced map from $\mathcal{M}_1^{*\mathfrak{w}_x}$ to X. Here is the definition: The map $\psi^{\mathfrak{w}_x}$ sends (C,p) from $\mathcal{M}_1^{*\mathfrak{w}_x}$ to the pair $(p, L = T_{1,0}C|_p)$. The first and third bullets in Proposition 3.5 imply that the differential of this map $\psi^{\mathfrak{w}_x}$ at $(C_x, x)$ is an isomorphism. (The differential is denoted subsequently as $D\psi^{\mathfrak{w}_x}$). The significance is that $\psi^{\mathfrak{w}_x}$ will be used to parametrize the submanifolds in $\mathcal{M}^{*\mathfrak{w}_x}$ near $C_x$. How this is done is explained momentarily.

To prepare for the explanation, note first that the tangent bundle to $\mathbb{P}(T_{1,0}X)$ is isomorphic to the direct sum $\mathfrak{p}^*TX \oplus \mathcal{V}$ with $\mathcal{V}$ being the kernel of the differential of the projection map to X (this is the map $\mathfrak{p}$). If $p \in X$ and C is a submanifold from $\mathcal{M}_{e,g}$ containing p, then this direct sum decomposition can be written at the point $(p, T_{1,0}C|_p)$ (which is a point in $\mathbb{P}(T_{1,0}X)$) as $TC|_p \oplus N_C|_p \oplus \mathcal{V}|_{(T_{1,0}C)|_p}$. The differential of $\psi^{\mathfrak{w}_x}$ at (C, p) can be viewed as linear map from $(T\mathcal{M}_1^{*\mathfrak{w}_x})|_{(C,p)}$ to this direct sum. To state the second point of note, let $\pi_{\mathcal{M}}$ denote the map from $\mathcal{M}_1^{*\mathfrak{w}_x}$ to $\mathcal{M}^{*\mathfrak{w}_x}$ sending any given (C, p) to C, and use $D\pi_{\mathcal{M}}$ to denote its differential. Note in this regard that the tangent space to $\mathcal{M}_1^{*\mathfrak{w}_x}$ at (C,p) is the space of pairs $(\eta, v)$ with $\eta$ from $\ker(D_{C,\mathfrak{w}})$ and $v \in TX|_p$ obeying $\eta|_p = \Pi_C v$. (Here again, $\Pi_C$ is the orthogonal projection to the normal bundle of C.) The map $D\pi_{\mathcal{M}}$ sends $(\eta, v)$ to $\eta$.

What follows is an algebraic consequence of what is said in the preceding paragraph: Supposing that $D\psi^{\mathfrak{w}_x}$ is an isomorphism at (C, p), the composition of first $(D\psi^{\mathfrak{w}_x})^{-1}$ and then $D\pi_{\mathcal{M}}$ restricts to the $N_C|_p \oplus \mathcal{V}|_{(T_{1,0}C)|_p}$ summand in $T\mathbb{P}(T_{1,0}X)|_p$ as an isomorphism from $N_C|_p \oplus \mathcal{V}|_{(T_{1,0}C)|_p}$ to $(T\mathcal{M}^{*\mathfrak{w}_x})|_C$. In the present context, Step 1 guarantees that $D\psi^{\mathfrak{w}_x}$ is an isomorphism at $(C_x, x)$. Therefore, the homomorphism $(D\pi_{\mathcal{M}})(D\psi^{\mathfrak{w}_x})^{-1}$ maps $N_{C_x}|_x \oplus \mathcal{V}|_{(T_{1,0}C_x)|_x}$ isomorphically to the tangent space to $\mathcal{M}^{*\mathfrak{w}_x}$ at $C_x$.

To explain the implications of the preceding observation, use the 1-forms $\sigma^0$ and $\sigma^1$ from (4.1) to define a product structure for $T^{1,0}X$ on $B_x$. This then identifies $\mathbb{P}(T_{1,0}X)$ on $B_x$ with the product of $B_x \times \mathbb{CP}^1$. Here is how: A point $\mathbb{CP}^1$ given by homogeneous coordinates $[\alpha, \beta]$ corresponds to the line in $T_{1,0}X$ at any given point in $B_x$ that is annihilated by the 1-form $-\alpha\sigma_0 + \beta\sigma_1$. For example, the tangent line to $C_x$ at $x$ corresponds to the point $[0, 1]$ in $\mathbb{CP}^1$. (In what follows, a disk about the origin in $\mathbb{C}$ is used to parametrize a disk about $[0, 1]$ in $\mathbb{CP}^1$ via the map that sends the complex coordinate $o$ of the disk to the point in $\mathbb{CP}^1$ with homogeneous coordinate $[o, 1]$. This then parametrizes the $\mathbb{C}$-valued 1-form $\sigma^1 - o\sigma^0$.)

Holding on to the preceding, introduce $\exp_{C_x}$ to denote the exponential map from a small radius disk bundle in $N_{C_x}$ to $X$ that is described in (2.13) of [T] or [T*]. Because $(D\pi_{\mathcal{M}})(D\psi^{\mathfrak{w}_x})^{-1}$ maps $N_{C_x}|_x \oplus \mathcal{V}|_{(T_{1,0}C_x)|_x}$ isomorphically to $(T\mathcal{M}^{*\mathfrak{w}_x})|_{C_x}$, the inverse function theorem provides the following:

- *A disk in $N_{C_x}|_x$ centered at $0$ (to be denoted by $D_x$ with coordinate $t$).*
- *A disk in $\mathbb{CP}^1$ centered at $[0, 1]$ (to be denoted by $\mathcal{O}_x$ and viewed as a disk about the origin in $\mathbb{C}$ with coordinate $o$).*
- *A diffeomorphism to be denoted by $\gamma_x$ from $D_x \times \mathcal{O}_x$ to an open neighborhood of $C_x$ in $\mathcal{M}^{*\mathfrak{w}_x}$ with compact closure taking $(0, 0)$ to $C_x$. This map takes any $(t, o) \in D_x \times \mathcal{O}_x$ to the submanifold in $\mathcal{M}^{*\mathfrak{w}_x}$ near $C_x$ that contains all entries of $\mathfrak{w}_x$ and the point $\exp_{C_x}(t)$, and whose tangent line at $\exp_{C_x}(t)$ is annihilated by the 1-form $\sigma^1 - o\sigma^0$.*
- *There exists $\kappa_x > 1$ such that if $\rho \in (0, \frac{1}{\kappa_x})$, and if the radii of both $D_x$ and $\mathcal{O}_x$ are less than $\rho$, then any submanifold in the image of $\gamma_x$ has $\eth$-distance less than $\kappa_x \rho$ from $C_x$.*

(4.2)

The notion of the $\eth$-distance between J-holomorphic subvarieties was introduce in (3.9) of [T] or [T*]: If $C$ and $C'$ are submanifolds, then $\eth(C; C') = 2\sup_{\{(z,z') \in C \times C'\}} \text{dist}(z, z')$.

Here is a schematic picture of $\gamma_x$:

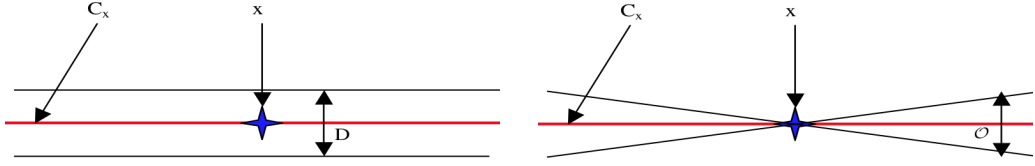

The left hand picture depicts submanifolds in $\mathcal{M}^{*w_x}$ obtained from $C_x$ by deformation via $\gamma_x$ in the $D_x$ direction. The right hand picture depicts submanifolds obtained by deformation in the $\mathcal{O}_x$ direction.

The submanifolds in the image of (4.2)'s map $\gamma_x$ can be depicted on $B_x$ with the help of the map $\lambda_x$ from (4.1) in the following way: The $\gamma_x$ image of any given point $(t = \tau, o = \hat{o})$ from $D_x \times \Delta_x$ appears in $B_x$ as the $\lambda_x$ image of the graph of a map from $\Delta_x$ to $D_x$ that can be written to leading order in $z$ and $\bar{z}$ as

$$z \to \tau + \hat{o}\, z - (\mathfrak{c}_t \tau + \mathfrak{c}_{\bar{t}} \bar{\tau})\bar{z} + \cdots$$

(4.3)

where $\mathfrak{c}_t$ and $\mathfrak{c}_{\bar{t}}$ are the respective $t$ and $\bar{t}$-derivatives at $t = 0$ of the function $\mathfrak{c}$ that appears in (4.1)'s the definition of the 1-from $\sigma^1$. The unwritten terms in (4.3) (which are indicated by $\cdots$) are bounded by $c_0 |z|^2 + \eth^2 |z|$ with $\eth$ bounded in turn by $c_0$ times the sum of the radii of $\Delta_x$ and $D_x$ and $\mathcal{O}_x$.

At points far from x, a picture of the image of $\gamma_x$ can be obtained by using a small radius ball in $\ker(D_{C_x})$ as coordinates on a neighborhood of $C_x$ in $\mathcal{M}_{e,g}$. What follows describes this larger coordinate neighborhood (this summarizes the discussion surrounding (2.14) in [T] or [T*]). Any given submanifold in $\mathcal{M}_{e,g}$ (call it C) has an open neighborhood in $\mathcal{M}_{e,g}$ consisting of submanifolds that are pointwise close C (their $\eth$-distance from C is small.) In particular, there exists such a neighborhood whose submanifolds are in 1-1 correspondence with submanifolds in X that be written as $\exp_C(\eta)$ with $\eta$ being a section of $N_C$ with small pointwise norm that obeys an equation that can be written as

$$D_C \eta + \mathfrak{r}_1 \cdot \partial \eta + \mathfrak{r}_0 = 0$$

(4.4)

with $\mathfrak{r}_1$ and $\mathfrak{r}_0$ as follows: These are smooth, fiber preserving maps from a small radius disk in the normal bundle of C to the respective bundles $\mathrm{Hom}(N \otimes T^{1,0}C; N \otimes T^{0,1}C)$ and to $N \otimes T^{0,1}C$ that obey $|\mathfrak{r}_1(b)| \leq c_0 |b|$ and $|\mathfrak{r}_0(b)| \leq c_0 |b|^2$. Meanwhile, the $L^2$-orthogonal projection map from $C^\infty(C; N_C)$ to the kernel of $D_C$ maps an open set of solutions to (4.4) centered on the zero section diffeomorphically to an open ball centered at 0 in $\ker(D_C)$. (In fact, the composition of this $L^2$ orthogonal projection map with the 1-1 correspondence to an open neighborhood of C in $\mathcal{M}_{e,g}$ is used to *define* the smooth structure on $\mathcal{M}_{e,g}$ on a neighborhood of C.) Supposing that C denotes a submanifold from $\mathcal{M}_{e,g}$, the convention employed henceforth (without further comment) uses the phrase '*an open neighborhood of C in $\mathcal{M}_{e,g}$*' to signify an open neighborhood of C that corresponds as just described to a small radius ball about the origin in $\ker(D_C)$.

With (4.2) in hand and the preceding understood, what follows is Step 2.

STEP 2: *Given* $x \in X$ *and given* $\mathfrak{w}_x$ *and* $C_x$ *from Step 1, fix disks* $D_x$ *and* $\mathcal{O}_x$ *as just described, with radius less than* $\frac{1}{\kappa_x}$ *(with* $\kappa_x$ *from (4.2)) so that the map* $\gamma_x$ *(with* $\gamma_x$ *as described by (4.2)) sends* $D_x \times \mathcal{O}_x$ *diffeomorphically onto a neighborhood of* $C_x$ *in* $\mathcal{M}^{*\mathfrak{w}_x}$.

### c) Step 3 of $\mathcal{K}$'s definition

The proposition that follows momentarily is referred to in the upcoming Step 3.

**Proposition 4.1**: *Given* $x \in X$ *and* $\mathfrak{w}_x$ *and* $C_x$ *from Step 1, there exists* $\kappa_x > 100$ *such that when* $\rho_x \in (0, \frac{1}{\kappa_x^2})$ *is used for the radii of the disks* $\Delta_x, D_x$ *and* $\mathcal{O}_x$, *then the maps* $\lambda_x$ *and* $\gamma_x$ *are defined using these disks, and all of what follows is true: Let* $\Delta_{x1}, D_{x1}$ *and* $\mathcal{O}_{x1}$ *denote the respective, concentric radius* $\frac{1}{\kappa_x} \rho_x$ *subdisks in* $\Delta_x, D_x$ *and* $\mathcal{O}_x$.

- *There is a submersion from* $\Delta_{x1} \times D_{x1} \times \mathcal{O}_{x1}$ *to* $D_x \times \mathcal{O}_x$ *to be denoted by* $\mathcal{F}_x$ *that is characterized as follows: Supposing that* $(z, \tau, \hat{o})$ *is from* $\Delta_{x1} \times D_{x1} \times \mathcal{O}_{x1}$, *then* $\gamma_x(\mathcal{F}_x(z, \tau, \hat{o}))$ *is the unique submanifold from the image of* $\gamma_x$ *that contains* $p = \lambda_x(z, \tau)$ *and whose* $(1, 0)$ *tangent line at* p *annihilates the 1-form* $\sigma^1 - \hat{o}\sigma^0$.
- *The map* $\mathcal{F}_x$ *differs from the canonical projection map* $(z, t, o) \to (t, o)$ *by at most* $\kappa_x(|z|^2 + |t|^2 + |o|^2)$. *Meanwhile, its differential differs from the differential of the projection by at most* $\kappa_x(|z| + |t| + |o|)$.
- *Fix* $(z, \tau, \hat{o})$ *from* $\Delta_{x1} \times D_{x1} \times \mathcal{O}_{x1}$ *again and let* $p = \lambda_x(z, \tau)$ *and let* $C_{p,\hat{o}} = \gamma_x(\mathcal{F}_x(z, \tau, \hat{o}))$.
  a) *The triple* $(C_{p,\hat{o}}, (p, \mathfrak{w}_x))$ *which is in* $\mathcal{M}_{e,g,d-1}$ *is a regular point of* $\pi_{d-1}$. *(This implies that* $(C_{p,\hat{o}}, \mathfrak{w}_x)$ *is a regular point of* $\pi_{d-2}$ *on* $\mathcal{M}_{e,g,d-2}$ *and that* $(C_{p,\hat{o}}, x)$ *is a regular point of* $\pi_1$ *on* $\mathcal{M}_{e,g,1}$.)
  b) *The submanifold* $C_{p,\hat{o}}$ *(which is in* $\mathcal{M}^{*(p,\mathfrak{w}_x)}$ *) is a regular point of the restriction of* $\psi^X$ *to* $\mathcal{M}^{*(p,\mathfrak{w}_x)}$.
  c) *Supposing that w is an entry of* $\mathfrak{w}_x$, *the submanifold* $C_{p,\hat{o}}$ *is a regular point for the map* $\psi^w$ *from* $\mathcal{M}^{*\mathfrak{w}_x}$ *to* $\mathbb{P}(T_{1,0}X)|_w$.
  d) *If* $q \in C_{p,\hat{o}}$ *is not an entry of* $\mathfrak{w}$, *then* $(C_{p,\hat{o}}, q)$ *is a regular point of* $\pi_1$ *on* $\mathcal{M}_1^{*\mathfrak{w}_x}$.

***Proof of Proposition 4.1***: The existence of $\mathcal{F}_x$ with the asserted properties follows from the inverse function theorem and (4.2). Meanwhile, the Items a)-d) of the third bullet follow because they hold for $(C_x, (x, \mathfrak{w}_x))$. Remember in this regard that Proposition 3.5 holds for $(C_x, (x, \mathfrak{w}_x))$, and because, as noted in Lemma 3.6, these items define open sets in the relevant moduli spaces.

What follows directly is a schematic picture of $C_{p,\hat{o}}$.

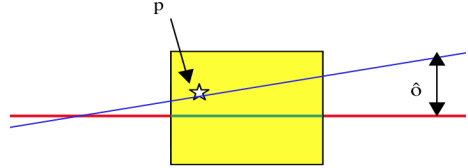

The blue line depicts $C_{p,\hat{o}}$. The double arrow depicts the angle ô.

To finish the set up for Step 3, for each $x \in X$, fix a number from $(0, \frac{1}{\kappa_x^2})$ (to be denoted by $\rho_x$) so that the conditions of Proposition 4.1 are satisfied with the radii of $\Delta_x$, $D_x$ (and thus $B_x$) and $\mathcal{O}_x$ equal to $\rho_x$. This number $\rho_x$ should be small enough so that $\lambda_x$ maps to a open neighborhood of $C_x$ whose submanifolds have distance at most $\kappa_x \rho_x$ from $C_x$. For each $x \in X$, let $B_{x4}$ denote the product of disks $\Delta_{x4} \times D_{x4}$ with $\Delta_{x4}$ denoting the concentric subdisk in $\Delta_x$ with radius $\frac{1}{4\kappa_x} \rho_x$, and with $D_{x4}$ denoting the concentric disk in $D_x$ with radius $\frac{1}{4\kappa_x} \rho_x$.

STEP 3: *Since X is compact, the cover $\{B_{x4}\}_{x \in X}$ has a finite subcover. Fix such a cover and let $\Xi \subset X$ denote the center points of the polydisks from this cover.*

Keep in mind that any given $x \in \Xi$ has its corresponding $\mathfrak{w}_x \in \mathcal{X}_{d-2}$ and $C_x \in \mathcal{M}_{e,g}$. The point x also has its versions of the disks $D_x \subset N_{C_x}|_x$ and $\mathcal{O}_x \subset \mathbb{P}(T_{1,0}X)|_x$; and it has the embedding $\gamma_x$ identifying $D_x \times \mathcal{O}_x$ with an open neighborhood in $\mathcal{M}^{*\mathfrak{w}_x}$ of the submanifold $C_x$. In addition, the point x has its corresponding map $\mathcal{F}_x$ from Proposition 4.1.

**d) Step 4 in $\mathcal{K}$'s definition**

Of concern now is the collection $\{\mathfrak{w}_x\}_{x \in \Xi}$. Since the conditions in Propositions and 4.1 are open conditions, the collection $\{\mathfrak{w}_x\}_{x \in \Xi}$ can be perturbed slightly (and the new set is henceforth denoted by the same labeling $\{\mathfrak{w}_x\}_{x \in \Xi}$) at the expense of shrinking the radii $\{\rho_x\}_{x \in \Xi}$ (slightly) that are used to define the polydisks that comprise the cover $\{B_x\}_{x \in \Xi}$. In particular, the collection $\{\mathfrak{w}_x\}_{x \in \Xi}$ and the radii $\{\rho_x\}_{x \in \Xi}$ can be chosen so that the following conditions are met:

- *Supposing that $x_1$, $x_2$ are distinct points from $\Xi$, then all entries of $\mathfrak{w}_{x_1}$ are distinct from all entries of $\mathfrak{w}_{x_2}$. No J-holomorphic subvariety in $\cup_{0 \leq k \leq g} \mathcal{M}_{e,k}$ contains all entries of both $\mathfrak{w}_{x_1}$ and $\mathfrak{w}_{x_2}$.*
- *If $x \in \Xi$, there is a submanifold $C_x \in \mathcal{M}_{e,g}$ containing x and each entry of $\mathfrak{w}_x$ as described in Step 1.*
- *If $x \in \Xi$, then the coordinate chart map $\lambda_x$ described by (4.1) is defined when $\Delta_x$ and $D_x$ have radius $\rho_x$.*
- *For each $x \in \Xi$, the diffeomorphism $\gamma_x$ from (4.2) is defined with $D_x$ and $\mathcal{O}_x$ have radius $\rho_x$; and the submanifolds in its image are in an open neighborhood of $C_x$ whose points have $\mathfrak{d}$-distance at most $\kappa_x \rho_x$ from $C_x$.*
- *The conclusions of each $x \in \Xi$ version of Proposition 4.1 holds with the given $\rho_x$.*
- *The set $\{B_{x4}\}_{x \in \Xi}$ defines an open cover of X. (Remember that if $x \in \Xi$, then $B_{x4}$ is the $\lambda_x$ image of $\Delta_{x4} \times D_{x4}$ with the latter being the product of the respective radius $\frac{1}{4\kappa_x}\rho_x$ concentric disks in $\Delta_x$ and $D_x$.*

(4.5)

The first bullet in (4.5) (which is the new item) follows from the first bullet of Proposition 3.7 in [T] or [T*]. Here is the fourth step

STEP 4: *Slightly move the points $\{\mathfrak{w}_x\}_{x \in \Xi}$ and (if necessary), shrink slighlty the radii $\{\rho_x\}_{x \in \Xi}$ in the manner just described so that all of the bullets in (4.5) hold.*

**e) Step 5 in $\mathcal{K}$'s definition**

The upcoming Step 5 also exploits the fact that the conditions in (4.5) are open conditions with regards to small changes in the chosen set $\{\mathfrak{w}_x\}_{x \in \Xi}$. To elaborate, it proves convenient to rename each $x \in \Xi$ version of Step 5's version of $\mathfrak{w}_x$ as $\mathfrak{w}_{x0}$. The essential point to be made is that the conditions in (4.5) hold when each $x \in \Xi$ version of $\mathfrak{w}_x$ is chosen from a small radius, open ball in $\mathcal{X}_{d-2}$ about the original choice $\mathfrak{w}_{x0}$. Furthermore, if the ball in $\mathcal{X}_{d-2}$ about any given $\mathfrak{w}_{x0}$ has small enough radius, then the disks $D_x$ and $\mathcal{O}_x$ and $\Delta_x$ can taken to be independent of $\mathfrak{w}_x$; and the $\mathfrak{w}_x$ version of the maps $\gamma_x$ and $\lambda_x$ will be very close (as an embedding of $D_x \times \mathcal{O}_x$ into $\mathcal{M}_{e,g}$) to the original $\mathfrak{w}_{x0}$ version of $\gamma_x$ and $\lambda_x$. In particular, there is an open ball in $\mathcal{X}_{d-2}$ centered on $\mathfrak{w}_{x0}$ (to be denoted by $\mathfrak{B}_x$) and an open neighborhood in $\mathcal{M}_{e,g}$ of the $\mathfrak{w}_{x0}$ version of $C_x$ with a submersion from $D_x \times \mathcal{O}_x \times \mathfrak{B}_x$ that is described by the next proposition. The notation is such that if $x \in \Xi$, then $\kappa_x$ is the number given by Proposition 4.1.

**Proposition 4.2**: *Fix* $x \in \Xi$, *and let* $\mathfrak{B}_x$ *denote a very small radius ball in* $X_{d-2}$ *centered at the point* $\mathfrak{w}_{x0}$. *There are smooth maps*

$$\Gamma_x: D_x \times \mathcal{O}_x \times \mathfrak{B}_x \to \mathcal{M}_{e,g} \quad \text{and} \quad \mathcal{P}_x: \Delta_{x1} \times D_{x1} \times \mathcal{O}_{x1} \times \mathfrak{B}_x \to D_x \times \mathcal{O}_x$$

*with the following properties:*
- *With regards to* $\Gamma_x$:
  a) $\Gamma_x(\cdot, \mathfrak{w}_{x0}) = \gamma_x(\cdot)$.
  b) *The map* $\Gamma_x$ *is a submersion onto an open neighborhood of* $C_x$ *in* $\mathcal{M}_{e,g}$ *with compact closure in* $\mathcal{M}_{e,g}$ *whose manifolds have* $\eth$-*distance at most* $\kappa_x \rho_x$ *from* C.
  c) *Supposing that* $(\tau, \hat{o}, \mathfrak{w})$ *is in* $D_x \times \mathcal{O}_x \times \mathfrak{D}_x$, *then* $\Gamma_x(\tau, \hat{o}, \mathfrak{w})$ *is the unique submanifold from the image of* $\Gamma_x$ *that contains each entry of* $\mathfrak{w}$ *and the point* $\exp_{C_x}(\tau)$, *and whose* $(1,0)$-*tangent space at the latter point is annihilated by the 1-form* $\sigma^1 - \hat{o}\sigma^0$.
- *With regards to* $\mathcal{P}_x$:
  a) $\mathcal{P}_x(\cdot, \mathfrak{w}_{x0}) = \mathcal{F}_x(\cdot)$.
  b) *Supposing that* $(z, \tau, \hat{o}, \mathfrak{w})$ *is from* $\Delta_{x1} \times D_{x1} \times \mathcal{O}_{x1} \times \mathfrak{B}_x$, *let* $p = \lambda_x(z, \tau)$ *and let* $C_{p,\hat{o},\mathfrak{w}}$ *denote* $\Gamma_x(\mathcal{P}(z,\tau,\hat{o},\mathfrak{w}), \mathfrak{w})$. *The submanifold* $C_{p,\hat{o},\mathfrak{w}}$ *is the unique submanifold from the image of* $\Gamma_x$ *that contains each entry of* $\mathfrak{w}$ *and the point* $p$, *and whose* $(1,0)$ *tangent line at* $p$ *is annihilated by the 1-form* $\sigma^1 - \hat{o}\sigma^0$.
  c) *The map* $\mathcal{P}_x$ *differs from the projection map (which send any given* $(z, t, o, \mathfrak{w})$ *to* $(t, o))$ *by at most* $\kappa_x(|z|^2 + |t|^2 + |o|^2 + |\mathfrak{w}|)$. *Meanwhile, its differential along the* $\Delta_{x1} \times D_{x1} \times \mathcal{O}_{x1}$ *differs from that of the projection by at most* $\kappa_x(|z| + |t| + |o| + |\mathfrak{w}|)$.
  d) *Fix* $(z, \tau, \hat{o}, \mathfrak{w})$ *from* $\Delta_{x1} \times D_{x1} \times \mathcal{O}_{x1} \times \mathfrak{B}_x$ *again and, again, let* $p = \lambda_x(z, \tau)$ *and let* $C_{p,\hat{o},\mathfrak{w}} = \Gamma_x(\mathcal{P}_x(z,\tau,\hat{o},\mathfrak{w}), \mathfrak{w})$.
     i) *The triple* $(C_{p,\hat{o},\mathfrak{w}}, p, \mathfrak{w})$ *which is in* $\mathcal{M}_{e,g,d-1}$ *is a regular point of* $\pi_{d-1}$.
     ii) *The submanifold* $C_{p,\hat{o},\mathfrak{w}}$ *which is in* $\mathcal{M}^{*(p,\mathfrak{w})}$ *is a regular point of the restriction of* $\psi^X$ *to* $\mathcal{M}^{*(p,\mathfrak{w})}$.
     iii) *Supposing that* $w$ *is an entry of* $\mathfrak{w}$, *the submanifold* $C_{p,\hat{o}}$ *is a regular point for the map* $\psi^w$ *to* $\mathbb{P}(T_{1,0}X)|_w$.
     iv) *If* $q \in C_{p,\hat{o},\mathfrak{w}}$ *is not an entry of* $\mathfrak{w}$, *then the pair* $(C_{p,\hat{o},\mathfrak{w}}, q)$ *is a regular point of* $\pi_1$ *on* $\mathcal{M}^{*\mathfrak{w}}_X$.

*Proof of Proposition 4.2*: All of the assertions follow directly from what is said in Proposition 4.1 using the inverse function theorem. In fact, the inverse function theorem says somewhat more about $\Gamma_x$. To elaborate, fix for the moment an entry, $w$, of $\mathfrak{w}_{x0}$. Let $\Delta_w \subset C_x$ denote a small radius disk centered at $w$ and let $D_w$ denote a small radius disk centered at 0 in the normal bundle to $C_x$ at $w$. With these defined, let $\lambda_w$ denote the coordinate chart map from $\Delta_w \times D_w$ to X that is described by (4.1) with x replaced by $w$. The polydisk $\Delta_w \times D_w$ is viewed using $\lambda_w$ as a small radius polydisk in X centered at $w$. Fix a point (to be denoted by $t_w$) in $\Delta_w$. Now let $\mathfrak{D}_x$ denote the product of the various versions of $\{t_w\} \times D_w$ as $w$ ranges through the d-2 entries of $\mathfrak{w}$. This product can be viewed using the various $\lambda_w$ coordinate chart maps as a polydisk in $\mathcal{X}_{d-2}$. Proposition 4.1 and the inverse function theorem implies that $\Gamma_x$ restricts to $D_x \times \mathcal{O}_x \times \mathfrak{D}_x$ as an embedding onto an open neighborhood of $C_x$ in $\mathcal{M}_{e,g}$.

Here is a parenthetical remark: The discussions that follows (and the discussion in Section 4 of [T] or [T*]) can make do with the replacement of $\mathfrak{B}_x$ in Proposition 4.2 and in (2.5) with $\mathfrak{D}_x$. Doing this requires only minor changes to what is done in Section 4 of [T] or [T*]. Even so, $\mathfrak{B}_x$ is used in what follows.

The following addendum to Proposition 4.2 is needed.

**Proposition 4.3**: *If all* $x \in \Xi$ *versions of Proposition 4.2's ball* $\mathfrak{B}_x$ *have small radius, then the following is true: Letting* $x_1, x_2$ *denote distinct points from* $\Xi$, *suppose that* $\mathfrak{w}_{x_1}$ *is from* $\mathfrak{B}_{x_1}$ *and that* $\mathfrak{w}_{x_2}$ *is from* $\mathfrak{B}_{x_2}$. *Then, no entry of* $\mathfrak{w}_{x_1}$ *is an entry of* $\mathfrak{w}_{x_2}$; *and no J-holomorphic subvariety in* $\cup_{0 \leq k \leq g} \mathcal{M}_{e,k}$ *contains all entries of both* $\mathfrak{w}_{x_1}$ *and* $\mathfrak{w}_{x_2}$.

*Proof of Proposition 4.3*: This follows from Proposition 3.7 in [T] or [T*] because the assertions hold for $\{\mathfrak{w}_{x0}\}_{x \in \Xi}$.

With the preceding lemmas in hand, what follows is Step 5.

STEP 5: *Fix a very small positive number to be denoted by s. Then, for each* $x \in \Xi$, *let* $\mathfrak{B}_x \subset \mathcal{X}_{d-2}$ *denote the radius s ball centered on* $\mathfrak{w}_{x0}$. *If s is sufficiently small (which is assumed henceforth), then the conclusions of Propositions 4.2 and 4.3 hold.*

**f) Step 6 in $\mathcal{K}$'s definition**

The final step that defines $\mathcal{K}$ follows directly. This step introduces by way of notation $D_{x2}$ and $\mathcal{O}_{x2}$ to denote the concentric disks in $D_x$ and $\mathcal{O}_x$ with radius $\frac{1}{2\kappa_x}\rho_x$.

STEP 6: *Having fixed a very small, positive number* $\rho$, *construct a finite set* $\Xi$ *as described in Step 3 with each* $x \in \Xi$ *version of* $\rho_x$ *being less than* $\rho$. *(The number* $\rho_x$ *for any given* $x \in \Xi$ *is chosen so that the conclusions of Proposition 4.1 hold.) Fix a very small, positive number s so that the conclusions of Propositions 4.2 and 4.3 hold when s* *is used to define the collection* $\{\mathfrak{B}_x\}_{x \in \Xi}$ *in Step 5. For each* $x \in \Xi$, *define the set* $\mathcal{K}_x$ *to be the* $\Gamma_x$-*image of* $D_{x2} \times \mathcal{O}_{x2} \times \mathfrak{B}_x$. *Having done this, set* $\mathcal{K} = \cup_{x \in \Xi} \mathcal{K}_x$.

An important point with regards to this definition is that there is a postive number, to be denoted by $r$, which is independent of $s$ (supposing that $s$ is small) and which is significant in the four ways that are listed below in Proposition 4.4. By way of a reminder with regards to notation from [T] or [T*]: Supposing that $x \in \Xi$ and that $\mathfrak{w} \in \mathfrak{B}_x$, then $\mathcal{Z}^{\mathfrak{w}}_X$ denotes the subset in $\mathcal{M}^{*\mathfrak{w}}_1$ of critical points of $\pi_1$ that map via $\pi_1$ to the complement in X of the entries of $\mathfrak{w}$. The $\pi_{\mathcal{M}}$ image of $\mathcal{Z}^{\mathfrak{w}}_X$ in $\mathcal{M}^{*\mathfrak{w}}$ is denoted by $\mathcal{Z}^{\mathfrak{w}}$. Supposing that $w$ is an entry of $\mathfrak{w}$, then $\mathcal{Y}^w$ denotes the set of critical points in $\mathcal{M}^{*\mathfrak{w}}$ of the map $\psi^w$. (Remember that this maps to $\mathbb{P}(T_{1,0}X)|_w$.) Let $\mathcal{Y}^{\mathfrak{w}}$ denote the union of these critical point loci indexed by the entries of $\mathfrak{w}$.

**Proposition 4.4**: *Given* $\Xi$, *and for each* $x \in \Xi$, *given the data* $(C_x, \mathfrak{w}_{x0}, \rho_x)$ *with* $\rho_x$ *small, there exists* $r > 0$ *which is such that if* $\mathcal{K}$ *is defined using a postive but very small value for* $s$, *then the following statements are true*:

- *If* x *and* x´ *are distinct points from* $\Xi$, *then each submanifold from* $\mathcal{K}_x$ *has* $\mathfrak{d}$-*distance greater than* $r$ *from each submanifold in* $\mathcal{K}_{x'}$.
- *Each submanifold in* $\mathcal{K}$ *has* $\mathfrak{d}$-*distance more than* $r$ *from any subvariety in* $\mathcal{M}_e - \mathcal{M}_{e,g}$.
- *If* $x \in \Xi$ *and if* $(\mathfrak{t}, \mathfrak{o}, \mathfrak{w}')$ *is from* $D_{x2} \times \mathcal{O}_{x2} \times \mathfrak{B}_x$, *then the submanifold* $\Gamma_x(\mathfrak{t}, \mathfrak{o}, \mathfrak{w}')$ *has* $\mathfrak{d}$-*distance greater than* $r$ *from any submanifold in any* $\mathfrak{w} \in \mathfrak{B}_x$ *version of* $\mathcal{Z}^{\mathfrak{w}}$.
- *If* $x \in \Xi$ *and if* $(\mathfrak{t}, \mathfrak{o}, \mathfrak{w}')$ *is from* $D_{x2} \times \mathcal{O}_{x2} \times \mathfrak{B}_x$, *then the submanifold* $\Gamma_x(\mathfrak{t}, \mathfrak{o}, \mathfrak{w}')$ *has* $\mathfrak{d}$-*distance greater than* $r$ *from any submanifold in any* $\mathfrak{w} \in \mathfrak{B}_x$ *version of* $\mathcal{Y}^{\mathfrak{w}}$.

*Proof of Proposition 4.4*: There is a choice for $r$ that makes the first bullet hold by virtue of what is said by Proposition 4.3. There is a choice for $r$ that makes the second bullet hold because $\Xi$ is finite and because each $x \in \Xi$ version of $\mathcal{K}_x$ has compact closure in $\mathcal{M}_{e,g}$. There is a choice for $r$ that makes the last two bullets hold because of what is said by Item d) of Proposition 4.2.

### 5. The current $\Phi$

Having chosen the numbers $r$ and $s$, define the current $\Phi$ as follows: Supposing that $\upsilon$ is a smooth 2-form on X, define $\Phi$ as follows:

$$\Phi(\upsilon) = \sum_{x \in \Xi} \int_{\{q \in D_{x2} \times \mathcal{O}_{x2} \times \mathcal{B}_x\}} \left( \int_{\Gamma_x(q)} \upsilon \right)$$

(5.1)

It is straightforward to verify that $\Phi$ defines a non-trivial current on X which is closed and non-negative and type 1-1. (This is to say that the current $\Phi$ defined above obeys the conclusions of Proposition 1.2 in [T] or [T*].) The follow proposition says more about $\Phi$ (it is the analog of Proposition 1.3 in [T] or [T*].)

**Proposition 5.1**: *Fix N >> 1 so as to define* $e = N[\omega]$. *If the almost complex structure is suitably generic, and if $\rho$ is small and then s is suitably small (given $\rho$), then the current $\Phi$ defined by (5.1) from the corresponding version of $\mathcal{K}$ has the following property: Fix t to be postive but small, and let $B \subset X$ denote a ball of radius t. Let $\sigma$ denote a unit length section of $T^{1,0}X|_B$ and let $\theta_B$ denote the characteristic function of B. Then*

$$\kappa^{-1} t^4 < \Phi(i \theta_B \sigma \wedge \bar\sigma) < \kappa t^4$$

*with $\kappa$ being positive, and independent of $\sigma$, B and t.*

*Proof of Proposition 5.1*: The proof uses much the same arguments as were used in the proof of Proposition 1.3 in [T] or [T*]. Because the detailed arguments for Proposition 5.1 are so similar to what is done in Section 4 of [T] or [T*], only a sketch is given below. There are two parts.

*Part 1*: This part considers the lower bound. To start the argument, take t so that the radius t ball about any given point in X is inside some $x \in \Xi$ version of $B_{x4}$ (which is the product of the concentric disks of radius $\frac{1}{4\kappa_x} \rho_x$ whose product is $B_x$). In particular, an upper bound for t can be chosen so that the distance from any point in any radius t ball to some $x \in \Xi$ version of $X - B_{x2}$ is at least 100t. With this understood, let B denote a given radius t ball and let $x \in \Xi$ denote a point whose version of $B_{x4}$ contains B in the manner just described.

Let p denote B's center point and let B´ denote the radius $\frac{1}{2}$ t ball with center p. Let $\mathcal{O}_{x2}$ denote the concentric, radius $\frac{1}{2\kappa_x} \rho_x$ disk in the disk $\mathcal{O}_x$. There is a point $o$ from $\mathcal{O}_{x2}$ with the following property: Let L denote the complex line in $\mathbb{P}(T_{1,0}X)|_p$ that is annihilated by $\sigma^1 - o\sigma^0$. Then

$$(i\sigma \wedge \bar\sigma)|_p(L) > c^{-1}$$

(5.2)

with $c$ here and in what follows denotes a number that is greater than 1 and independent of the radius t of B, the center p of B, and the 1-form $\sigma$. The value of $c$ can be assumed to increase between sucessive appearances. Since (5.2) holds for the center point of B, then it follows (by continuity) that (5.2) holds with p replaced by any point in B (if $t < c^{-1}$) and when L is annihilated by $\sigma^1 - \hat{o}\sigma^0$ with $\hat{o}$ being any point from the radius $c^{-1}$ disk centered at the point $o$. Denote this disk by $\hat{O}$.

Reintroduce the maps $\Gamma_x$ and $\mathcal{P}_x$ from Proposition 4.2. Write p using $\lambda_x$ as $\lambda_x(z, \tau)$ with $(z, \tau)$ from $B_{x2}$. Supposing that $\mathfrak{w} \in \mathfrak{B}_x$, then the submanifold $\Gamma_x(\mathcal{P}_x(z,\tau,\hat{o},\mathfrak{w}),\mathfrak{w})$ is the unique submanifold from $\mathcal{K}_x$ that contains p, all entries of $\mathfrak{w}$ and whose tangent line at p annihilates $\sigma^1 - \hat{o}\sigma^0$. And, if $\xi \in \mathbb{C}$ has norm at most $\frac{1}{100}t$, then $\Gamma_x(\mathcal{P}_x(z,\tau+\xi,\hat{o},\mathfrak{w}),\mathfrak{w})$ is a submanifold from $\mathcal{K}_x$ that also intersects the radius $\frac{1}{2}t$ ball centered at p. Let C denote any one of these submanifolds with $|\xi| < \frac{1}{100}t$ and with $\hat{o}$ being a point in the disk $\hat{O}$. According to Lemma 2.2 in [T] or [T*], and by virtue of (5.2),

$$\int_C \theta_B i\sigma \wedge \bar{\sigma} \geq c^{-1}t^2 .$$

(5.3)

Meanwhile, it is a consequence of what is said by Proposition 4.2 that the subspace inside $D_{x2} \times \mathcal{O}_{x2} \times \mathfrak{B}_x$ that parametrizes this $(\xi, \hat{o}, \mathfrak{w})$ family of submanifolds via $\Gamma_x$ and $\mathcal{P}_x$ has volume greater than $c^{-1}t^2 s^{4d-8}$. To explain: A factor uniformly proportional to $t^2$ come from the area of the disk of possible choices for $\xi$; a factor uniformly bounded away from zero comes from the area of the disk $\hat{O}$; and a factor uniformly proportional to $s^{4d-8}$ comes from the volume of $\mathfrak{B}_x$.

This volume weight $c^{-1}t^2 s^{4d-8}$ times the $c^{-1}t^2$ lower bound in (5.3) leads directly to a $c^{-1}t^4 s^{4d-8}$ lower bound for $\Phi(i\theta_B \sigma \wedge \bar{\sigma})$.

*Part 2*: Consider next the upper bound assertion in Proposition 5.1. The arguments for this upper bound differ very little from what is said in Sections 4.5 and 4.6 of [T] (or Sections 4e and 4f of [T*]). Note in this regard that an upper bound for the number $\Phi(i\theta_B \sigma \wedge \bar{\sigma})$ proportional to $t^4$ (and otherwise independent of B, and independent of $\sigma$) holds if there exists $c_* > 1$ (which should be independent of t, B and $\sigma$) such that if $x \in \Xi$, and supposing that $s < c_*^{-1}$ has been fixed and that $t < c_*^{-1}s$, then:

$$\int_{\{q \in D_{x2} \times \mathcal{O}_{x2} \times \mathfrak{B}_x\}} (\int_{\Gamma_x(q)} \theta_B i\sigma \wedge \bar{\sigma}) \leq c_* t^4 .$$

(5.4)

The four steps that follow sketch the argument for the existence of $c_*$.

<u>Step 1</u>: The proof that there exists such a number $c_*$ requires first a fact from Lemma 2.2 in [T] or [T*]: The area in B of any submanifold from $\mathcal{M}_{e,g}$ is bounded by $ct^2$. This implies directly that

$$\int_C \theta_B i\sigma \wedge \bar{\sigma} \leq ct^2 .$$

(5.5)

when $|\sigma| \leq 1$. Therefore, the bound in (5.4) follows from (5.5) with a proof of that there exists $c_\Diamond > 1$ which is independent of t, B and $\sigma$ such that

*The set of* $q \in D_{x2} \times \mathcal{O}_{x2} \times \mathfrak{B}_x$ *with* $\Gamma_x(q) \cap B \neq \emptyset$ *has measure at most* $c_\Diamond t^2$ .

(5.6)

The rest of the steps sketch the argument from Section 4 of [T] or [T*] for the existence of the desired $c_\Diamond$ when s is small.

<u>Step 2</u>: The proof that there is a $c_\Diamond$ for (5.6) exploits the fact that any given submanifold in the image of $\Gamma_x$ can be written as $\exp_{C_x}(\eta)$ with $\eta$ being a small normed section of the normal bundle to $C_x$ that obeys (4.4). Supposing that $q = (t, o, \mathfrak{w})$ is in $D_{x2} \times \mathcal{O}_{x2} \times \mathfrak{B}_x$, then the corresponding version of $\eta$ is written as $\eta_q$. This implies (among other things) that if B intersects a submanifold from the image of $\Gamma_x$, then the whole of B necessarily lies in the $\exp_{C_x}$ embedded image of a small radius disk bundle about the zero section in the normal bundle to $C_x$.

<u>Step 3</u>: Consider first the case when the center of B has distance greater than 4s from any entry of $\mathfrak{w}_{x0}$. With this restriction understood, fix $\mathfrak{w} \in \mathfrak{B}_x$ and introduce by way of notation $C_{B,\mathfrak{w}}$ to denote the subset of $D_{x2} \times \mathcal{O}_{x2}$ that is mapped by $\Gamma_x(\cdot, \mathfrak{w})$ to a submanifold that intersects B. The bound in (5.6) in this case is a consequence of:

*If* $s < c^{-1}$ *and if* $\mathfrak{w} \in \mathfrak{B}_x$, *then the volume of* $C_{B,\mathfrak{w}}$ *is at most* $ct^2$ .

(5.7)

To prove (5.7), let p now denote the center point of B. This point lies in the $\exp_{C_x}$ image of a disk centered at 0 in a unique fiber of the normal bundle to $C_x$. Let $\pi(p)$ denote the base point in C of this fiber disk. A map to be denoted by $f$ from $D_{x1} \times \mathcal{O}_{x1}$ to $N_{C_x}|_{\pi(p)}$ is defined as follows: Supposing that $(t, o) \in D_{x1} \times \mathcal{O}_{x1}$, and setting $q = (t, o, \mathfrak{w})$, take $f(t, o)$ to be the point $\eta_q(\pi(p))$. Let $Df$ denote the differential of $f$.

The fourth bullet of Proposition 4.4 implies that the map $f$ is a submersion (this is the raison-d'etre of this bullet). Therefore, because $\pi(p)$ is uniformly far from any entry of $\mathfrak{w}$ (distance at least $4s$), the norm of $Df$ on the orthogonal complement of kernel($Df$) is bounded below by $c^{-1}$. This has the following consequence (by appeal to the implicit function theorem): If $D'$ is a small radius disk in $N_{C_x}|_{\pi(p)}$, then the volume in $D_{x2} \times \mathcal{O}_{x2}$ of $f^{-1}(D')$ is bounded by $c$ times the area of $D'$. This is relevant to the problem at hand because if $t < c^{-1}$, then there is a disk in $N_{C_x}|_{\pi(p)}$ of radius at most $ct^2$ that contains the intersection point of every submanifold from $C_{B,\mathfrak{w}}$ (these are the submanifolds from the $\Gamma_x(\cdot, \mathfrak{w})$ image of $D_{x2} \times \mathcal{O}_{x2}$ that hit B). Therefore, the volume of $C_{B,\mathfrak{w}}$ is bounded by $ct^2$.

<u>Step 4</u>: Consider next the case when B's center point has distance $4s$ or less from some entry of $\mathfrak{w}_{x0}$. Nothing is lost by assuming that this is first entry of $\mathfrak{w}$. It is denoted in what follows by $w_0$. (Note that if $s$ is small, then there is only one such entry.) The function giving the distance to B's center point defines a function on the ball of radius $8s$ centered on $w_0$. Denote this function by $\ell$. If $w$ is in the radius $s$ ball centered at $w_0$ and if $\mathfrak{w} = (w, \mathfrak{w}')$ is a point in $X \times X_{d-3}$ from $\mathfrak{B}_x$, then the area in $\mathbb{P}(T_{1,0}X)|_w$ of submanifolds from $\Gamma_x(\cdot, \mathfrak{w})$ that intersect B is no greater than $c \frac{t^2}{\ell^2}$ (see (3.6)). It follows from this, and from the fourth bullet of Proposition 4.4 that the volume of the set of submanifolds in the $\Gamma_x(\cdot, \mathfrak{w})$ image of $D_{x2} \times \mathcal{O}_{x2}$ that intersect B is also bounded by $c \frac{t^2}{\ell^2}$. (This appeal is the reason for the fourth bullet of Proposition 4.4.) Therefore, the volume of the set of $\mathfrak{q} \in \mathfrak{B}$ with $\Gamma_x(\mathfrak{q})$ intersecting B is no greater than $c$ times the integral of $\frac{t^2}{\ell^2}$ over the ball of radius $s$ centered at $w$. This, in turn, is no greater than

$$c \int_0^{4s} \frac{t^2}{\ell^2} \ell^3 \, d\ell$$

(5.8)

which is bounded by $ct^2$ as required.

**References**

[G]    M. Gromov, *Pseudoholomorphic curves in symplectic manifolds*, Invent. Math. **82** (1985), 307–347.

[LL1]  T.-J. Li and A. Liu, *General wall crossing formula*. Math. Res. Lett. **2** (1995) 797-810.

[LL2]  T.-J. Li and A. Liu, *Uniqueness of the symplectic class, surface cone and symplectic cone of 4-manifolds with $b^{2+} = 1$*. J. Differential Geom. **58** (2001) 331-370.

[LL3]  T-J. Li and A. Liu, *The equivalence between SW and Gr in the case where $b^+ = 1$*. Internat. Math. Res. Notices **7** (1999) 335–345.

[T]  C. H. Taubes, *Tamed to compatible: Symplectic forms via moduli space integration*. J. Symplectic Geom. **9** (2011) 161-250. See also arXiv:0910.5440v2 for a recent version with minor corrections. This arXiv version is denoted by reference [T*] in the preceding text. Note that the designation of subsections is not the same in the arXiv version and the published version.